\documentclass[12pt]{article}

\usepackage{amsfonts,amsmath}
\usepackage{latexsym}

\author{ K\'aroly J. B\"or\"oczky\footnote{Supported by
NKFIH grants K 132002}, Apratim De}

\title{Stability of the Pr\'ekopa-Leindler inequality for log-concave functions}

\newcommand{\proof}{\noindent{\it Proof: }}
\newcommand{\proofbox}{\mbox{ $\Box$}\\}
\newcommand{\R}{\mathbb{R}}

\newtheorem{lemma}{LEMMA}[section]
\newtheorem{theo}[lemma]{THEOREM}

\newtheorem{coro}[lemma]{COROLLARY}

\newtheorem{prop}[lemma]{PROPOSITION}

\begin{document}

\maketitle

\begin{abstract}
A stability version of the
 Pr\'ekopa-Leindler inequality for log-concave functions on $\R^n$ is established.
\end{abstract}

\noindent MSC Subject Index: 26D15

\section{Introduction}

For $X\subset \R^n$, we write ${\rm conv}\,X$ to denote the convex hull of $X$, and say that $X$ is homothetic to 
$Y\subset\R^n$  if $Y=\gamma X+z$ for $\gamma>0$ and $z\in\R^n$.
Writing $|X|$ to denote Lebesgue measure of a measurable subset $X$ of $\R^n$
(with $|\emptyset|=0$), the Brunn-Minkowski inequality 
(Schneider \cite{Sch14}) says that if $\alpha,\beta>0$ and $X,Y,Z$ are bounded measurable subsets of $\R^n$, then
\begin{equation}
\label{BrunnMinkowski}
|Z|^{\frac1n}\geq \alpha |X|^{\frac1n}+\beta|Y|^{\frac1n}
\mbox{  \ provided \ } \alpha X+\beta Y \subset Z,
\end{equation}
and in the case $|X|,|Y|> 0$, equality holds if and only if
${\rm conv}\,X$ and ${\rm conv}\,Y$  are homothetic convex bodies with
$|({\rm conv}\,X)\backslash X|=|({\rm conv}\,Y)\backslash Y|=0$ and
${\rm conv}\,Z=\alpha({\rm conv}\,X)+\beta ({\rm conv}\,Y)$. We note that 
even if $X$ and $Y$ are Lebesgue measurable, the Minkowski linear combination 
$\alpha X+\beta Y$ may not be measurable.

Because of the homogeneity of the Lebesgue measure, an equivalent form of \eqref{BrunnMinkowski} is the following. 
If $\lambda\in(0,1)$, then
\begin{equation}
\label{BrunnMinkowskiprod}
|Z|\geq |X|^{1-\lambda}|Y|^{\lambda}
\mbox{  \ provided \ } (1-\lambda)X+\lambda Y \subset Z.
\end{equation}
In the case $|X|,|Y|> 0$,  equality in \eqref{BrunnMinkowskiprod} implies that
${\rm conv}\,X$ and ${\rm conv}\,Y$  are translates, and $|({\rm conv}\,X)\backslash X|=|({\rm conv}\,Y)\backslash Y|=0$.

For convex $X$ and $Y$, the first stability forms of the Brunn-Minkowski inequality were due to Minkowski himself 
(see Groemer \cite{Gro93}).
If the distance of the convex $X$ and $Y$ is measured in terms of the so-called Hausdorff distance, then
Diskant \cite{Dis73} and Groemer \cite{Gro88} provided close to optimal stability versions
(see Groemer \cite{Gro93}). However, the natural distance is in terms of the volume of the symmetric difference, and the optimal result is due to Figalli, Maggi, Pratelli \cite{FMP09,FMP10}.
To define the ``homothetic distance'' $A(K,C)$
of convex bodies $K$ and $C$, let $\alpha=|K|^{\frac{-1}n}$ and
$\beta=|C|^{\frac{-1}n}$, and let
$$
A(K,C)=\min\left\{|\alpha K\Delta (x+\beta C)|:\,x\in\R^n\right\}.
$$
In addition, let
$$
\sigma(K,C)=\max\left\{\frac{|C|}{|K|},\frac{|K|}{|C|}\right\}.
$$

\begin{theo}[Figalli, Maggi, Pratelli]
\label{Maggi}
For $\gamma^*(n)=(\frac{(2-2^{\frac{n-1}{n}})^{\frac32}}{122n^7})^2$,
and  for any convex bodies $K$ and $C$ in $\R^n$,
$$
|K+C|^{\frac1n}\geq (|K|^{\frac1n}+|C|^{\frac1n})
\left[1+\frac{\gamma^*}{\sigma(K,C)^{\frac1n}}\cdot A(K,C)^2\right].
$$
\end{theo}

Here the exponent $2$ of $A(K,C)^2$ is optimal, see Figalli, Maggi, Pratelli \cite{FMP10}.
We note that prior to \cite{FMP10}, the only known error term in the Brunn-Minkowski inequality
was of order $A(K,C)^\eta$ with $\eta\geq n$,
due to 
Diskant \cite{Dis73} and  Groemer \cite{Gro88} 
(see Groemer \cite{Gro93}). For a more direct approach, we refer to
Esposito, Fusco, Trombetti \cite{EFT05}.

The Figalli, Maggi, Pratelli \cite{FMP09} factor of the form $\gamma^*(n)=cn^{-14}$ for some absolute constant $c>0$ was improved to $c n^{-7}$ by Segal \cite{Seg12}, and subsequently to $c n^{-5.5}$ by Kolesnikov, Milman \cite{KoM}, Theorem~12.12. 
The current best known bound for $\gamma^*(n)$ is $n^{-5 - o(1)}$, which follows by combining the general estimate of Kolesnikov-Milman \cite{KoM}, Theorem 12.2, with the bound $n^{o(1)}$ on the Cheeger constant of a convex body in isotropic position that follows from Yuansi Chen’s work \cite{Che21} on the Kannan-Lovasz-Simonovits conjecture.
Harutyunyan \cite{Har18} conjectured that $\gamma^*(n)=c n^{-2}$ is the optimal order of the constant, and 
showed that it can't be of smaller order. Actually, Segal \cite{Seg12} proved that
 Dar's conjecture in \cite{Dar99} would imply that we may choose $\gamma^*(n)=c n^{-2}$
for some absolute constant $c>0$. 

The paper Eldan, Klartag \cite{ElK14} discusses "isomorphic" stability versions of the Brunn-Minkowski inequality under condition of the type
$|\frac12\,K+\frac12\,C|\leq 5\sqrt{|K|\cdot|C|}$,
and considers, for example, the $L^2$ Wasserstein distance of the uniform measures on suitable affine images of $K$ and $C$.

If $X$ is measurable bounded and $Y$ is a convex body, then improving on the estimate in Carlen, Maggi \cite{CaM17},
the paper  Barchiesi, Julin \cite{BaJ17} proves that
\begin{equation}
\label{XconvY}
|X+Y|^{\frac1n}\geq |X|^{\frac1n}+|Y|^{\frac1n}+\delta_n\min\{|X|,|Y|\}^{\frac1n}A(X,Y)^2
\end{equation}
for some $\delta_n>0$ depending on $n$.

Let us shortly discuss the case when $X,Y,Z$ are bounded measurable with positive measure and
$X+Y\subset Z$ (and hence neither $X$ nor $Y$ are assumed to be convex). In this case,
estimates similar to \eqref{XconvY} (having $Z$ in place of $X+Y$) were proved
by Hintum, Spink, Tiba \cite{HST} if $X=Y$ and $n\geq 1$, and  by Hintum, Spink, Tiba \cite{HST2} 
if $n=2$ and $X,Y$ any bounded measurable set (and even better error term of order $A(X,Y)$ if $n=1$  is proved by Freiman, see Christ \cite{Chr12}). If $n\geq 3$ and $X,Y,Z$ any bounded measurable sets with $X+Y\subset Z$, then only a much weaker estimate is known;
namely, Figalli, Jerison \cite{FiJ15,FiJ17} prove that
if
$$
|\,|X|-1|+|\,|Y|-1|+|\,|Z|-1|<\varepsilon
\mbox{  \ and \ $\frac12\,X+\frac12\,Y \subset Z$}
$$
for small $\varepsilon>0$,
then there exist a convex body (i.e. compact convex sets with non-empty interior)  $K$ and $z\in \R^n$ such that
$$
\mbox{$X\subset K$, \ $Y+z\subset K$ \ and \ $|K\backslash X|+|K\backslash (Y+z)|<c_n\varepsilon^\eta$}
$$
where $c_n,\eta>0$ depend on $n$ and $\eta<n^{-3^n}$.

Our main theme is a generalization of the Brunn-Minkowski inequality known as Pr\'ekopa-Leindler inequality. 
The inequality itself, due to
Pr\'ekopa \cite{Pre71} and Leindler \cite{Lei72} in dimension one, was generalized in
 Pr\'ekopa \cite{Pre73,Pre75} and  Borell \cite{Bor75} ({\it cf.} also Marsiglietti \cite{Mar17}, Bueno, Pivovarov \cite{BuP}), Brascamp, Lieb \cite{BrL76}, Kolesnikov, Werner \cite{KoW}, 
Bobkov, Colesanti, Fragal\`a \cite{BCF14}).
Various applications are provided and surveyed in Ball \cite{Bal},
Barthe \cite{Bar}, 
Fradelizi, Meyer \cite{FrM07} and Gardner \cite{Gar02}.
The following multiplicative version from \cite{Bal} is often more useful
and is more convenient for geometric applications.

\begin{theo}[Pr\'ekopa-Leindler]
\label{PLn}
If $\lambda\in(0,1)$ and $h,f,g$ are non-negative integrable functions on $\R^n$
satisfying $h((1-\lambda)x+\lambda y)\geq f(x)^{1-\lambda}g(y)^\lambda$ for
$x,y\in\mathbb{R}^n$, then
\begin{equation}
\label{PLnineq}
\int_{\R^n}  h\geq \left(\int_{\R^n}f\right)^{1-\lambda} \cdot \left(\int_{\R^n}g\right)^\lambda.
\end{equation}
\end{theo}

It follows from Theorem~\ref{PLn} that the Pr\'ekopa-Leindler inequality has the following multifunctional form
which resembles Barthe's Reverse Brascamp-Lieb inequality \cite{Bar98}.
If $\lambda_1,\ldots,\lambda_m>0$ satisfy $\sum_{i=1}^m\lambda_i=1$ and 
$f_1,\ldots,f_m$ are non-negative integrable functions on $\R^n$, then
\begin{equation}
\label{PLnmulti}
\int_{\R^n}^*
\sup_{z=\sum_{i=1}^m\lambda_ix_i}
  \prod_{i=1}^mf(x_i)^{\lambda_i}\,dz
\geq  \prod_{i=1}^m \left(\int_{\R^n}f_i\right)^{\lambda_i}
\end{equation}
where $*$ stands for outer integral in the case the integrand is not measurable.

We say that a function $f:\,\R^n\to[0,\infty)$ have {\it positive integral} if $f$ is measurable and
$0<\int_{\R^n}f<\infty$. For a convex subset $\Gamma\subset \R^n$, we say that a 
function $f:\,\Gamma\to[0,\infty)$ is {\it log-concave}, if for any $x,y\in\Gamma$ and
$\alpha,\beta\in[0,1]$ with $\alpha+\beta=1$, we have
$f(\alpha x+\beta y)\geq f(x)^\alpha g(y)^\beta$.
The case of equality in Theorem~\ref{PLn} has been characterized by
 Dubuc \cite{Dub77}.

\begin{theo}[Dubuc]
\label{PLnequa}
If $\lambda\in(0,1)$ and $h,f,g:\,\R^n\to[0,\infty)$ have positive integral,
satisfy $h((1-\lambda)x+\lambda y)\geq f(x)^{1-\lambda}g(y)^\lambda$ for
$x,y\in\mathbb{R}^n$ and equality holds in \eqref{PLnineq}, then
$f,g,h$ are log-concave up to a set of measure zero, and
there exist $a>0$ and $z\in \R^n$ such that 
\begin{eqnarray*}
f(x)&=& a^\lambda\,h(x-\lambda z)\\
g(x)&=& a^{-(1-\lambda)}h(x+(1-\lambda) z)
\end{eqnarray*}
for almost all $x$.
\end{theo}

Our goal is to prove a stability version of the Pr\'ekopa-Leindler inequality Theorem~\ref{PLn} at least for log-concave functions.

\begin{theo}
\label{PLhstab0}
For some absolute constant $c>1$, if $\tau\in(0,\frac12]$, $\tau\leq \lambda\leq 1-\tau$, 
$h,f,g:\,\R^n\to [0,\infty)$ are integrable  such that
   $h((1-\lambda)x+\lambda\,y)\geq f(x)^{1-\lambda}g(y)^\lambda$ for
$x,y\in\mathbb{R}^n$, $h$ is log-concave and
$$
\int_{\R^n}  h\leq (1+\varepsilon) \left(\int_{\R^n}f\right)^{1-\lambda}\left( \int_{\R^n}g\right)^\lambda
$$
for $\varepsilon\in(0,1]$, then
there exists $w\in\R^n$ such that setting $a=\int_{\R^n}g/\int_{\R^n}f$, we have
\begin{eqnarray*}
\int_{\R^n}|f(x)-a^{-\lambda} h(x-\lambda\,w)|\,dx&\leq &c^nn^n\sqrt[19]{\frac{\varepsilon}{\tau}}  
\cdot \int_{\R^n}f \\
\int_{\R^n}|g(x)-a^{1-\lambda}h(x+(1-\lambda)w)|\,dx&\leq &c^nn^n\sqrt[19]{\frac{\varepsilon}{\tau}}   
\cdot\int_{\R^n}g.
\end{eqnarray*}
\end{theo}
\noindent{\bf Remark } According to Lemma~\ref{log-concave-geometric-mean} (i), if $f$ and $g$ are log-concave, then
$$
h(z)=\sup_{z=(1-\lambda)x+\lambda y}f(x)^{1-\lambda}g(y)^\lambda
$$
is log-concave, as well, and hence Theorem~\ref{PLhstab0} applies. \\

A statement similar to Theorem~\ref{PLhstab0} was proved by Ball, B\"or\"oczky \cite{BB11} in the case of even log-concave functions and $\tau=\frac12$ with error term of order $\varepsilon^{\frac16}|\log \varepsilon|^{\frac23}$ instead of 
the $\varepsilon^{\frac1{19}}$ in Theorem~\ref{PLhstab0}. We note that 
Bucur, Fragal\`a \cite{BuF14} proved a nice stability version of the Pr\'ekopa-Leindler inequality for log-concave functions
if the distance of functions is measured not in terms of the "translative" $L_1$ distance, but the weaker notion of bounding the (translative) distance of all one dimensional projections.

An "isomorphic" stability result for the Prekopa-Leindler inequality, in terms of the transportation distance is obtained in Eldan \cite{Eld13}, Lemma~5.2. By rather standard considerations, one can show that non-isomorphic stability results in terms of transportation distance imply stability in terms of $L_1$ distance (e.g., such implication is attained by combining Proposition~2.9 in Bubeck, Eldan, Lehec \cite{BEL18} and Proposition~10 in Eldan, Klartag \cite{ElK14}). However, the current result in \cite{Eld13}, due to its isomorphic nature, falls short of being able to obtain a meaningful bound in terms of the $L_1$ distance.

Brascamp, Lieb \cite{BrL76} proved a local version of the Prekopa-Leindler inequality for log-concave functions (Theorem~4.2 in \cite{BrL76}), which is equivalent to a Poincare-type so called Brascamp-Lieb inequality Theorem~4.1 in \cite{BrL76}. 
The paper Livshyts \cite{Liv} provides a stability version of this Brascamp-Lieb inequality.

Let us present a version of Theorem~\ref{PLhstab0} analogous to
Theorem~\ref{Maggi}. If $f,g$ are non-negative functions on $\R^n$ with
$0<\int_{\R^n}f<\infty$ and $0<\int_{\R^n}g<\infty$, then for the probability densities
$$
\tilde{f}=\frac{f}{\int_{\R^n}f} \mbox{ \ and \ }\tilde{g}=\frac{g}{\int_{\R^n}g},
$$
we define 
\begin{equation}
\label{tildeL1}
\widetilde{L}_1(f,g)=\inf_{v\in\R^n}\int_{\R^n}|\tilde{f}(x-v)-\tilde{g}(x)|\,dx.
\end{equation}

\begin{coro}
\label{PLhstablikeMaggi}
If  $\tau\in(0,\frac12]$,
$\lambda\in[\tau,1-\tau]$ and
$f,g$ are log-concave functions with positive integral on $\R^n$, then
$$
\int_{\R^n}\sup_{z=(1-\lambda)x+\lambda y} f(x)^{1-\lambda}g(y)^\lambda\,dz
\geq \left(1+\gamma\cdot \tau\cdot\widetilde{L}_1(f,g)^{19}\right) \left(\int_{\R^n}f\right)^{1-\lambda}\left( \cdot \int_{\R^n}g\right)^\lambda
$$
where $\gamma=c^n/n^{19n}$ for some absolutute constant $c\in(0,1)$.
\end{coro}

We also deduce a stability version of \eqref{PLnmulti} from Theorem~\ref{PLhstab0}.

\begin{theo}
\label{PLhstabcor0}
For some absolute constant $c>1$, if  $\tau\in(0,\frac1m]$, $m\geq 2$, 
$\lambda_1,\ldots,\lambda_m\in[\tau,1-\tau]$ satisfy $\sum_{i=1}^m\lambda_i=1$ and 
$f_1,\ldots,f_m$ are log-concave functions with positive integral on $\R^n$ such that
$$
\int_{\R^n}\sup_{z=\sum_{i=1}^m\lambda_ix_i} \prod_{i=1}^mf_i(x_i)^{\lambda_i}\,dz
\leq (1+\varepsilon)  \prod_{i=1}^m \left(\int_{\R^n}f_i\right)^{\lambda_i}
$$
for $\varepsilon\in(0,1]$, then for the log-concave 
$h(z)=\sup_{z=\sum_{i=1}^m\lambda_ix_i} \prod_{i=1}^mf(x_i)^{\lambda_i}$,
there exist $a_1,\ldots,a_m>0$ and  $w_1,\ldots,w_m\in\R^n$ such that
$\sum_{i=1}^m\lambda_iw_i=o$ and
for $i=1,\ldots,m$, we have
$$
\int_{\R^n}|f_i(x)-a_ih(x+w_i)|\,dx\leq c^nn^nm^5
\sqrt[19]{\frac{\varepsilon}{m\tau}}\cdot
\int_{\R^n}f_i.
$$
\end{theo}
\noindent{\bf Remark } $a_i=\frac{\left(\int_{\R^n}f_i\right)^{1-\lambda_i}}
{\prod_{j\neq i}\left(\int_{\R^n}f_j\right)^{\lambda_j}}$ for $i=1,\ldots,m$
in Theorem~\ref{PLhstabcor}.\\

See Corollary~\ref{log-concave-geometric-mean-cor} for the log-concavity of the $h$ in 
Theorem~\ref{PLhstabcor}.

Recently, various breakthrough stability results about geometric functional inequalities have been obtained. 
Fusco, Maggi, Pratelli \cite{FMP08} proved an optimal stability version of the isoperimetric inequality (whose result was extended to the Brunn-Minkowski inequality by Figalli, Maggi, Pratelli \cite{FMP09,FMP10}). 
 Stonger versions of the 
functional Blaschke-Santal\'o inequality is provided by 
Barthe, B\"or\"oczky, Fradelizi \cite{BBF14}, of the Borell-Brascamp-Lieb inequality is provided by Ghilli, Salani \cite{GhS17}, Rossi, Salani \cite{RoS17,RoS19} and Balogh, Krist\'aly \cite{BaK18} (later even on Riemannian manifolds),
 of the Sobolev inequality by Figalli, Zhang \cite{FiZ} (extending Bianchi, Egnell \cite{BiE91} and Figalli, Neumayer \cite{FiN19}), 
Nguyen \cite{Ngu16} and Wang \cite{Wan16}, of the log-Sobolev inequality by Gozlan \cite{Goz},
and of some related inequalities by 
Caglar, Werner \cite{CaW17}, Cordero-Erausquin \cite{Cor17}, Kolesnikov, Kosov \cite{KoK17}.

The next Section~\ref{seconevariable} reviews the known stability versions of the Prekopa-Leindler inequality for functions on $\R$, and  Section~\ref{secidea} outlines the idea of the proofs of 
Theorem~\ref{PLhstab0}, Corollary~\ref{PLhstablikeMaggi} and Theorem~\ref{PLhstabcor0}.

\section{Stability versions of the one dimensional Pr\'ekopa-Leindler inequality}
\label{seconevariable}

If $n=1$, then Ball, B\"or\"oczky \cite{BB10}
  provided the following stability version of the Prekopa-Leindler inequality Theorem~\ref{PLn}
in the logconcave case.

\begin{theo}
\label{PLstab}
There exists a positive absolute constant $c$ with
the following property:
If $h,f,g$ are non-negative integrable functions
with positive integrals on $\R$ such that
$h$ is log-concave,  $h(\frac{r+s}2)\geq \sqrt{f(r)g(s)}$ for
$r,s\in\mathbb{R}$,  and
$$
\int_{\R}  h\leq (1+\varepsilon) \sqrt{\int_{\R}f \cdot \int_{\R}g},
$$
for $\varepsilon\in(0,1)$, then
there exists $b\in\R$ such that for $a=\sqrt{\int_{\R}g/\int_{\R}f}$, we have
\begin{eqnarray*}
\int_{\R}|f(t)-a\,h(t+b)|\,dt&\leq &
c\cdot\sqrt[3]{\varepsilon}|\ln \varepsilon|^{\frac43}\cdot
\int_{\R}f(t)\,dt \\
\int_{\R}|g(t)-a^{-1}h(t-b)|\,dt&\leq &
c\cdot\sqrt[3]{\varepsilon}|\ln \varepsilon|^{\frac43}\cdot
\int_{\R}g(t)\,dt.
\end{eqnarray*}
\end{theo}
{\bf Remark } If $f$ and $g$ are log-concave probability distributions,
then $a=1$, and if in addition $f$ and $g$ have
the same expectation,
then even $b=0$ can be assumed.\\

We note that combining Theorem~\ref{PLstab}, Lemma~\ref{log-concave-geometric-mean} (ii) and 
Lemma~\ref{log-concave-error} implies the following more precise stability version of the one-dimensional Prekopa-Leindler inequality.

\begin{coro}
\label{PLstabdim1lambda}
For some absolute constant $c>1$, if $\tau\in(0,\frac12]$, $\tau\leq \lambda\leq 1-\tau$, 
$h,f,g:\,\R\to [0,\infty)$ are integrable  such that
   $h((1-\lambda)x+\lambda\,y)\geq f(x)^{1-\lambda}g(y)^\lambda$ for
$x,y\in\mathbb{R}^n$, $h$ is log-concave and
$$
\int_{\R^n}  h\leq (1+\varepsilon) \left(\int_{\R^n}f\right)^{1-\lambda}\left( \int_{\R^n}g\right)^\lambda
$$
for $\varepsilon\in(0,1]$, then
there exists $w\in\R$ such that for $a=\int_{\R^n}g/\int_{\R^n}f$, we have
\begin{eqnarray*}
\int_{\R^n}|f(x)-a^{-\lambda} h(x-\lambda\,w)|\,dx&\leq &
c\left(\frac{\varepsilon}{\tau}\right)^{\frac13} |\log \varepsilon|^{\frac43} 
\cdot \int_{\R^n}f \\
\int_{\R^n}|g(x)-a^{1-\lambda}h(x+(1-\lambda)w)|\,dx&\leq &
c\left(\frac{\varepsilon}{\tau}\right)^{\frac13} |\log \varepsilon|^{\frac43} 
\cdot\int_{\R^n}g.
\end{eqnarray*}
\end{coro}

As it was observed by C. Borell \cite{Bor75},
and later independently by K.M. Ball \cite{Bal},
assigning to any function $H:[0,\infty]\to[0,\infty]$
the function  $h:\R\to[0,\infty]$ defined by
$h(x)=H(e^x)e^x$,
we have the version Theorem~\ref{PLB} of the Pr\'ekopa-Leindler
inequality. We note that if $H$
is log-concave and decreasing, then $h$
 is log-concave.

\begin{theo}
\label{PLB}
If $H,F,G:[0,\infty]\to[0,\infty]$  integrable functions
satisfy  $H(\sqrt{rs})\geq \sqrt{F(r)G(s)}$ for
$r,s\geq 0$, then
$$
\int_0^\infty  H\geq \sqrt{\int_0^\infty F\cdot \int_0^\infty G}.
$$
\end{theo}

Therefore we deduce the following statement by
Theorem~\ref{PLstab}:

\begin{coro}
\label{PLBstab}
There exists a positive absolute constant $c_0>1$ with
the following property:
If $H,F,G:[0,\infty]\to[0,\infty]$ are integrable functions
with positive integrals such that
$H$ is log-concave and decreasing,  $H(\sqrt{rs})\geq \sqrt{F(r)G(s)}$ for
$r,s\in[0,\infty]$,  and
$$
\int_0^\infty  H\leq (1+\varepsilon)
\sqrt{\int_0^\infty F \cdot \int_0^\infty G}
$$
for $\varepsilon\in[0,c_0^{-1})$, then
there exist $a,b>0$, such that
\begin{eqnarray*}
\int_0^\infty|F(t)-a\,H(b\,t)|\,dt&\leq &
c\cdot\sqrt[3]{\varepsilon}|\ln \varepsilon|^{\frac43}\cdot
\int_0^\infty F(t)\,dt \\
\int_0^\infty|G(t)-a^{-1}H(b^{-1}t)|\,dt&\leq &
c\cdot\sqrt[3]{\varepsilon}|\ln \varepsilon|^{\frac43}\cdot
\int_0^\infty G(t)\,dt.
\end{eqnarray*}
\end{coro}
{\bf Remark }
If in adddition,
$F$ and $G$ are decreasing log-concave probability distributions
then $a=b$ can be assumed. The condition that
$H$ is log-concave and decreasing can be replaced by the
one that $H(e^t)$ is log-concave.\\

Concerning general measurable functions, at least the stability of the one-dimensional Brunn-Minkowski inequality has been clarified by Christ \cite{Chr12}
(see also Theorem~1.1 in Figalli, Jerison \cite{FiJ17}).

\begin{theo}
\label{dim1BMstab}
If $X,Y\subset \R$ are measurable with $|X|,|Y|>0$, and $|X+Y|\leq |X|+|Y|+\delta$ for some  
$\delta\leq \min\{|X|,|Y|\}$,  then there exist intervals $I,J\subset\R$ such that
$X\subset I$, $Y\subset J$, $|I\backslash X|\leq \delta$ and $|J\backslash Y|\leq \delta$.
\end{theo}

\section{Ideas to verify Theorem~\ref{PLhstab0} and its consequences}
\label{secidea}

For Theorem~\ref{PLhstab0}, the main  goal is to prove Theorem~\ref{PLhstab12} which is essentially the case $\lambda=\frac12$ of Theorem~\ref{PLhstab0}
for log-concave functions and for small $\varepsilon$,
and then the general case is handled in Sections~\ref{seclambda} and \ref{secThCor}.

\begin{theo}
\label{PLhstab12}
If $h,f,g:\,\R^n\to [0,\infty)$ are log-concave, $f,g$ are probability distributions,  $h(\frac{x+y}2)\geq \sqrt{f(x)g(y)}$ for
$x,y\in\mathbb{R}^n$,  and
$$
\int_{\R^n}  h\leq 1+\varepsilon
$$
where $0<\varepsilon<(cn)^{-n}$, then
there exists $w\in\R^n$ such that
\begin{eqnarray*}
\int_{\R^n}|f(x)-h(x-w)|\,dx&\leq& \tilde{c}n^8\cdot \sqrt[18]{\varepsilon}\cdot|\log\varepsilon|^{n}\\
\int_{\R^n}|g(x)-h(x+w)|\,dx&\leq &\tilde{c} n^8\cdot\sqrt[18]{\varepsilon}\cdot|\log\varepsilon|^{n} 
\end{eqnarray*}
where $c,\tilde{c}>1$ are absolute constants.
\end{theo}

Our proof of the stability version Theorem~\ref{PLhstab12} of the Pr\'ekopa-Leindler inequality stems from the Ball's following
argument ({\it cf}  \cite{Bal} and Borell \cite{Bor75}) proving the Pr\'ekopa-Leindler inequality based on the Brunn-Minkowski inequality. 

Let $f,g,h:\R^n\to[0,\infty]$ have positive integrals and satisfy that
$h(\frac{x+y}2)\geq \sqrt{f(x)g(y)}$ for
$x,y\in\mathbb{R}^n$, and for $t>0$, let
\begin{eqnarray*}
\Phi_t&=&
\{x\in\R^n:\,f(x)\geq t\} \mbox{ \ and \ } F(t)=|\Phi_t|\\
\Psi_t&=&\{x\in\R^n:\,g(x)\geq t\} \mbox{ \ and \ } G(t)=|\Psi_t|\\
\Omega_t&=&\{x\in\R^n:\,h(x)\geq t\}\mbox{ \ and \ } H(t)=|\Omega_t|.
\end{eqnarray*}
As it was observed by Ball \cite{Bal} and and Borell \cite{Bor75}, the condition
on $f,g,h$ yields that if  $\Phi_r,\Psi_s\neq\emptyset$ for $r,s>0$, then
$$
\mbox{$\frac12$}(\Phi_r+\Psi_s)\subset \Omega_{\sqrt{rs}}.
$$
Therefore the Brunn-Minkowski inequality implies that
$$
H(\sqrt{rs})\geq \left(\frac{F(r)^{\frac1n}+G(s)^{\frac1n}}2\right)^n
\geq\sqrt{F(r)\cdot G(s)}
$$
for all $r,s>0$.
In particular we deduce the Pr\'ekopa-Leindler inequality
by Theorem~\ref{PLB}, as
$$
\int_{\R^n}h=\int_0^\infty H(t)\,dt\geq
\sqrt{\int_0^\infty F(t)\,dt \cdot \int_0^\infty G(t)\,dt}
=\sqrt{\int_{\R^n}f \cdot \int_{\R^n}g}.
$$

Turning to Theorem~\ref{PLhstab12}, we will need the stability version Corollary~\ref{PLBstab}
of the Prekopa-Leindler inequality for function on $\R$, 
and a stability version of the product form of the Brunn-Minkowski inequality on $\R^n$.
Since
\begin{eqnarray*}
\frac12\left(|K|^{\frac1n}+|C|^{\frac1n}\right)&=&
|K|^{\frac1{2n}}|C|^{\frac1{2n}}
\left[1+\frac12
\left(\sigma(K,C)^{\frac1{4n}}-\sigma(K,C)^{\frac{-1}{4n}}\right)^2\right] \\
&\geq &|K|^{\frac1{2n}}|C|^{\frac1{2n}}
\left[1+\frac{(\sigma(K,C)-1)^2}{32n^2\sigma(K,C)^{\frac{4n-1}{2n}}}\right],
\end{eqnarray*}
using the notation $\sigma=\sigma(K,C)=\max\{\frac{|C|}{|K|},\frac{|K|}{|C|}\}$,
we conclude from the stability version Theorem~\ref{Maggi} of the Brunn-Minkowski inequality by Figalli, Maggi, Pratelli \cite{FMP10} that
\begin{equation}
\label{Maggiprod}
\left|\mbox{$\frac12$}(K+C)\right|\geq \sqrt{|K|\cdot|C|}\left[1+
\frac{(\sigma-1)^2}{32n\sigma^2}
+\frac{n\gamma^*(n)}{\sigma^{\frac1n}}\cdot A(K,C)^2\right].
\end{equation}

We observe that  the volume of the symmetric difference $|K\Delta C|$ of
 convex bodies $K$ and $C$ is a metric on convex bodies in $\R^n$.
We use this fact in the following consequence of Theorem~\ref{Maggi}:

\begin{lemma}
\label{VKC1L}
If $\eta\in(0,\frac1{122n^7})$ and $K,C.L$ are convex bodies in $\R^n$ such that
$|C|=|K|$, $|L|\leq (1+\eta)|K|$  and
$\frac12\,K+\frac12\,C\subset L$, then there exists $w\in \R^n$ such that
$$
|K\Delta(L-w)|\leq 245n^7\sqrt{\eta}\,|K|\mbox{ \ and \ }
|C\Delta(L+w)|\leq 245n^7\sqrt{\eta}\,|K|.
$$ 
\end{lemma}
\proof We may assume that $|C|=|K|=1$.
According to Theorem~\ref{Maggi}, there exists $z\in\R^n$, such that  
$$
|K\cap (C-z)|\geq 1-\sqrt{\frac{\eta}{\gamma^*(n)}}>1-122n^7\sqrt{\eta}.
$$
It follows from $z+[K\cap (C-z)]\subset C$ that 
$M=\frac12\,z+[K\cap (C-z)]\subset L$, and hence $|L|\leq 1+\eta$ implies
$|L\Delta M|<\eta+122n^7\sqrt{\eta}<123n^7\sqrt{\eta}$.
Writing $w=\frac12\,z$, we have
$$
|K\Delta(L-w)|\leq |K\Delta(M-w)|+|(M-w)\Delta(L-w)|< 245n^7\sqrt{\eta},
$$
 and similar argument yields $|C\Delta (L+w)|< 245n^7\sqrt{\eta}$.
\proofbox

To  prove Theorem~\ref{PLhstab12}, first we  discuss some fundamental estimates for log-concave functions in
Section~\ref{seclog-concave}, then compare the level sets of $f$, $g$  and $h$ in Theorem~\ref{PLhstab12}
in Section~\ref{seclevelsets}, and finally complete the argument for Theorem~\ref{PLhstab12} in Section~\ref{sechalf}.
We verify Theorem~\ref{PLhstab0} when $\varepsilon$ is small in Section~\ref{seclambda},
and prove Theorem~\ref{PLhstab0} and Corollary~\ref{PLhstablikeMaggi}
in Section~\ref{secThCor}. Finally, Theorem~\ref{PLhstabcor0} is proved in Section~\ref{secmanyfunctions}.

\section{Some properties of log-concave functions}
\label{seclog-concave}

First, we characterize a log-concave function $\varphi$ on $\R^n$, $n\geq 2$ with positive integral; namely, if
$0<\int_{\R^n}\varphi<\infty$. For any measurable function $\varphi$ on $\R^n$, we define
$$
M_\varphi=\sup \varphi.
$$

\begin{lemma}
\label{log-concave-positive-integral}
Let $\varphi:\,\R^n\to[0,\infty)$ be log-concave. Then $\varphi$ has positive integral if and only if 
$\varphi$ is bounded, $M_\varphi>0$, and for any $t\in(0,M_f)$, the level set 
$\{\varphi>t\}$ is bounded and has non-empty interior.
\end{lemma}
\proof If $\varphi$ has positive integral, then $M_\varphi>0$, and there exists some $t_0\in(0,M_\varphi)$ such that the $n$-dimensional measure of $\{\varphi>t_0\}$ is positive. As $\{\varphi>t_0\}$ is convex, it has non-empty interior. It follows from the log-concavity of $\varphi$ that the level set $\{\varphi>t\}$ has non-empty interior for any 
$\varphi\in(0,M_f)$. In turn, we deduce that $\varphi$ is bounded from  the log-concavity of $\varphi$ and $\int_{\R^n}\varphi<\infty$.

Next we suppose that that there exists $t\in(0,M_\varphi)$ such that the level set 
$\{\varphi>t\}$ is unbounded and seek a contradiction. As $\{\varphi>t\}$ is convex, there exists a 
$u\in S^{n-1}$ such that $x+su\in {\rm int}\{\varphi>t\}$
for any $x\in {\rm int}\{\varphi>t\}$ and $s\geq 0$. We conclude that $\int_{\R^n}\varphi=\infty$, contradicting the assumption
$\int_{\R^n}\varphi<\infty$. Therefore the level set $\{\varphi>t\}$ is bounded for any $t\in(0,M_\varphi)$.

Assuming that the conditions of Lemma~\ref{log-concave-positive-integral} hold, we readily have
$\int_{\R^n}\varphi>0$. To show $\int_{\R^n}\varphi<\infty$, we choose $x_0\in\R^n$ such that $\varphi(x_0)>0$, and let
$B$ be an $n$-dimensional ball of centered $x_0$ and radius $\varrho>0$ containing 
$\{\varphi>\frac1e\, \varphi(x_0)\}$. Let us consider
$$
\psi(x)=\varphi(x_0)e^{-\frac{\|x-x_0\|}\varrho}.
$$
It follows from the log-concavity of $\varphi$ that $\varphi(x)\leq \psi(x)$ if $\|x-x_0\|\geq \varrho$, and hence
$$
\int_{\R^n}\varphi\leq \int_{B}\varphi+\int_{\R^n\backslash B}\psi< \infty,
$$
verifying Lemma~\ref{log-concave-positive-integral}.
\proofbox

For a measurable bounded function $\varphi$ on $\R^n$ and for $t\in\R$, let
$$
\Xi_{\varphi,t}=\{x\in\R^n:\,\varphi(x)\geq t\}.
$$
If $\varphi$ is  log-concave with positive integral, then
we consider the symmetric decreasing rearrangement
$\varphi^*:\,\R^n\to \R$ where
$$
|\Xi_{\varphi,t}|=|\Xi_{\varphi^*,t}|
$$
for any $t>0$, and if $|\Xi_{\varphi,t}|>0$, then $\Xi_{\varphi^*,t}$ is a Euclidean ball centered at the origin $o$,
and
$$
M_\varphi=\max_{x\in\R^n} \varphi(x)=\max_{x\in\R^n} \varphi^*(x)=\varphi^*(o).
$$
It follows from Lemma~\ref{log-concave-positive-integral} that $\varphi^*$ is well-defined.
We deduce that $\varphi^*$ is also a log-concave, and
$$
\int_{\R^n}\varphi=\int_0^\infty |\Xi_{\varphi,t}|\,dt=\int_0^\infty |\Xi_{\varphi^*,t}|\,dt=\int_{\R^n}\varphi^*.
$$

We write $B^n$ to denote the unit Euclidean ball in $\R^n$ centered at the origin, and $\kappa_n=|B^n|$, and hence
the surface area of $S^{n-1}$ is $n\kappa_n$. 
For log-concave functions, a useful property of the  symmetric decreasing rearrangement is
that if $\varrho\,B^n=\Xi_{\varphi^*,s\,M_\varphi}$ 
and $s=e^{-\gamma\,\varrho}$ for $\gamma,\varrho>0$, then
\begin{equation}
\label{phi*z}
\begin{array}{rcl}
\varphi^*(x)&\geq &M_\varphi e^{-\gamma\|x\|}\mbox{ \ provided $\|x\|\leq \varrho$}\\
\varphi^*(x)&\leq &M_\varphi e^{-\gamma\|x\|}\mbox{ \ provided $\|x\|\geq \varrho$}.
\end{array}
\end{equation}
For $s=e^{-\gamma\,\varrho}$, we have
\begin{equation}
\label{Xisz}
|\Xi_{\varphi,sM_\varphi}|=|\Xi_{\varphi^*,sM_\varphi}|=\kappa_n \varrho^n.
\end{equation}

As a related integral, it follows from induction on $n$ that
\begin{equation}
\label{0infty}
\int_0^{\infty}e^{-\gamma r}r^{n-1}\,dr=(n-1)!\cdot \gamma^{-n}.
\end{equation}

\begin{lemma}
If $\varphi$ is a log-concave probability density on $\R^n$, then 
\begin{equation}
\label{closetomax}
|\Xi_{\varphi,(1-\tau)M_\varphi}|\geq \mbox{$\frac1{n!+1}\,$}\tau^nM_\varphi^{-1}\mbox{ \ 
for $\tau\in(0,1)$}. 
\end{equation}
\end{lemma}
\proof
To prove (\ref{closetomax}) based on (\ref{phi*z}), let $\gamma,\varrho>0$ be such that
$\varrho\,B^n=\Xi_{\varphi^*,(1-\tau)\,M_\varphi}$ 
and $1-\tau=e^{-\gamma\,\varrho}$. 
Here $e^{-\gamma\,\varrho}>1-\gamma\,\varrho$ yields $\gamma\,\varrho\geq \tau$, thus
it follows from (\ref{Xisz}) and  (\ref{0infty})  that
\begin{eqnarray*}
1&= &\int_{\R^n}\varphi^*(x)\,dx\leq
|\varrho\,B^n|\cdot M_\varphi+
\int_{\R^n\backslash\varrho\,B^n}\varphi^*(x)\,dx\\
&\leq &|\Xi_{\varphi,(1-\tau)M_\varphi}|\cdot M_\varphi+
\int_{\R^n}M_\varphi e^{-\gamma\,\|x\|}\,dx\\
&=&M_\varphi\cdot |\Xi_{\varphi,(1-\tau)M_\varphi}|+
M_\varphi n\kappa_n\int_0^{\infty}e^{-\gamma r}r^{n-1}\,dr\\
&=&M_\varphi\cdot |\Xi_{\varphi,(1-\tau)M_\varphi}|+
M_\varphi n!\kappa_n\cdot \gamma^{-n}\\
&\leq & M_\varphi\cdot |\Xi_{\varphi,(1-\tau)M_\varphi}|+
M_\varphi n!\kappa_n\cdot \frac{\varrho^n}{\tau^n}\\
&=&M_\varphi\cdot |\Xi_{\varphi,(1-\tau)M_\varphi}|\left(1+\frac{n!}{\tau^n}\right),
\end{eqnarray*}
proving (\ref{closetomax}).
\proofbox

We note that the estimate in \eqref{closetomax} is close to be optimal because if
the probability density is of the form
$\varphi(x)=M_f e^{-\gamma\|x\|}$ for suitable $\gamma>0$, then
$|\Xi_{\varphi,(1-\tau)M_\varphi}|= \frac{|\ln (1-\tau)|^n}{n!M_\varphi}<\frac{e\tau^n}{n!M_\varphi}$
if $\tau\in(0,\frac1n)$.

For a log-concave probability density $\varphi$, let $\mu_\varphi$ be the probability measure associated to 
$\varphi$; namely, $d\mu_\varphi(x)=\varphi(x)\,dx$.
According to Lemma~5.16 in Lov\'asz, Vempala \cite{LoV07}, if $s\in(0,e^{-4(n-1)})$, then 
\begin{equation}
\label{muphiss}
\mu_\varphi(\varphi<sM_\varphi)\leq \frac{e^{n-1}}{(n-1)^{n-1}}\cdot s\cdot |\ln s|^{n-1}\leq
s\cdot |\ln s|^n.
\end{equation}
But what we really need is the following estimate.

\begin{lemma}
If $s\in(0,e^{-4(n-1)})$ and $\varphi$ is a log-concave probability density on $\R^n$, then
\begin{eqnarray}
\label{philessthans0}
|\Xi_{\varphi,sM_\varphi}|&< & 
\frac{2|\ln s|^n}{n!M_\varphi},\\
\label{philessthans}
\int_0^{sM_\varphi}|\Xi_{\varphi,t}|\,dt&<&\left(1+\frac{1}{M_\varphi}\right) s\cdot |\ln s|^n.
\end{eqnarray}
\end{lemma}
\proof To prove (\ref{philessthans}) based on (\ref{phi*z}), let $\gamma,\varrho>0$ be such that
$\varrho\,B^n=\Xi_{\varphi^*,s\,M_\varphi}$ 
and $s=e^{-\gamma\,\varrho}$.
Since $s\in(0,e^{-4(n-1)})$, we deduce that
\begin{equation}
\label{gammavarrhon-1}
\gamma\varrho>4(n-1).
\end{equation}
It follows from \eqref{0infty} and integration by parts that
\begin{equation}
\label{intfromrho}
\int_{\varrho}^{\infty}e^{-\gamma r}r^{n-1}\,dr=
e^{-\gamma\varrho}\int_{0}^{\infty}e^{-\gamma r}r^{n-1}\,dr\cdot 
\sum_{k=0}^{n-1}\frac{(\gamma\varrho)^k}{k!}.
\end{equation}
Here the well-known estimate $(n-1)!>\frac{(n-1)^{n-1}}{e^{n-1}}$ implies that if $k=1,\ldots,n-1$, then
\begin{equation}
\label{n-1n-k}
(n-1)\cdot\ldots\cdot (n-k)>\frac{(n-1)^k}{e^k}.
\end{equation}
In addition, 
$1+e\cdot s<e^{\frac34\,s}$ holds for $s\geq 4$.
Combining this with \eqref{n-1n-k} yields that if $t>4(n-1)$, then 
\begin{equation}
\label{sumt0ton-1}
\sum_{k=0}^{n-1}\frac{t^k}{k!}< \sum_{k=0}^{n-1}{n-1\choose k}\left(\frac{et}{n-1}\right)^k=
\left(1+\frac{et}{n-1}\right)^{n-1}<e^{\frac34\,t}.
\end{equation}
Therefore, we deduce from \eqref{gammavarrhon-1}, \eqref{intfromrho} and \eqref{sumt0ton-1} that
\begin{eqnarray}
\nonumber
\int_{\varrho}^{\infty}e^{-\gamma r}r^{n-1}\,dr&<&
e^{-\gamma\varrho}\int_{0}^{\infty}e^{-\gamma r}r^{n-1}\,dr\cdot
e^{\frac34\,\gamma\varrho}=e^{-\frac14\,\gamma\varrho}\int_{0}^{\infty}e^{-\gamma r}r^{n-1}\,dr\\
\label{rhoinfty}
&<&
\frac12\int_{0}^{\infty}e^{-\gamma r}r^{n-1}\,dr.
\end{eqnarray}
In particular, using \eqref{phi*z}, \eqref{rhoinfty} and later \eqref{Xisz}, we have
\begin{eqnarray*}
1&\geq &\int_{\varrho\,B^n}\varphi^*(x)\,dx\geq \int_{\varrho\,B^n}
M_\varphi e^{-\gamma\|x\|}\,dx\\
&=&n\kappa_nM_\varphi\int_0^{\varrho}e^{-\gamma r}r^{n-1}\,dr\\
&\geq &n\kappa_nM_\varphi\cdot \frac12\,\int_0^{\infty}e^{-\gamma r}r^{n-1}\,dr=
\frac{n!M_\varphi\kappa_n\varrho^n}{2\gamma^n\varrho^n}\\
&=&\frac{n!M_\varphi}{2}\cdot 
\frac{|\Xi_{\varphi,s}|}{|\ln s|^n}.
\end{eqnarray*}
We conlude that if $s\in(0,e^{-4(n-1)})$, then
$$
|\Xi_{\varphi,s}|\leq 
\frac{2|\ln s|^n}{n!M_\varphi}.
$$
Combining the last inequality with (\ref{muphiss}), we deduce that if $s\in(0,e^{-4(n-1)})$, then
$$
\int_0^{sM_\varphi}|\Xi_{\varphi,t}|\,dt=
|\Xi_{\varphi,s}|\cdot s+\mu_\varphi(\varphi<sM_\varphi) <
e\left(1+\frac{1}{M_\varphi}\right) s\cdot |\ln s|^n,
$$
proving \eqref{philessthans}. \proofbox

We note that the estimate in \eqref{philessthans0} is close to be optimal because if
again the probability density is of the form
$\varphi(x)=M_f e^{-\gamma\|x\|}$ for suitable $\gamma>0$, then
$|\Xi_{\varphi,sM_\varphi}|= \frac{|\ln s|^n}{n!M_\varphi}$.

\section{The area of the level sets in Theorem~\ref{PLhstab12}}
\label{seclevelsets}

Let $f,g,h$ as in Theorem~\ref{PLhstab12}.
We may assume 
that 
\begin{equation}
\label{fogo}
f(o)=\max\{f(x):x\in\R^n\}\mbox{ and }g(o)=\max\{g(x):x\in\R^n\}.
\end{equation}
According to  Lemma~\ref{log-concave-positive-integral}, for $t>0$, we may consider the bounded convex sets
\begin{eqnarray*}
\Phi_t&=&\Xi_{f,t}=
\{x\in\R^n:\,f(x)\geq t\} \mbox{ \ and \ } F(t)=|\Phi_t|\\
\Psi_t&=&\Xi_{g,t}=\{x\in\R^n:\,g(x)\geq t\} \mbox{ \ and \ } G(t)=|\Psi_t|\\
\Omega_t&=&\Xi_{h,t}=\{x\in\R^n:\,h(x)\geq t\}\mbox{ \ and \ } H(t)=|\Omega_t|
\end{eqnarray*}
where  \eqref{fogo} yields
\begin{equation}
\label{oinPhitPsit}
o\in \Phi_t\cap \Psi_t,
\end{equation}
and we have
\begin{equation}
\label{FGint1}
\int_0^\infty F=\int_0^{M_f} F=\int_{\R^n} f=1\mbox{ \ and \ }
\int_0^\infty G=\int_0^{M_g} G=\int_{\R^n} g=1.
\end{equation}

As it was observed in K.M. Ball \cite{Bal}, the condition
on $f,g,h$ yields that if  $\Phi_r,\Psi_s\neq\emptyset$ for $r,s>0$, then
\begin{equation}
\label{sections}
\mbox{$\frac12$}(\Phi_r+\Psi_s)\subset \Omega_{\sqrt{rs}}.
\end{equation}
Therefore the Brunn-Minkowski inequality yields that
\begin{equation}
\label{minksum}
H(\sqrt{rs})\geq \left(\frac{F(r)^{\frac1n}+G(s)^{\frac1n}}2\right)^n
\geq\sqrt{F(r)\cdot G(s)}
\end{equation}
for all $r,s>0$.

Let $c_0>1$ be the absolute constant of Corollary~\ref{PLBstab}, and if $0<\varepsilon<1/c_0$, then let
\begin{equation}
\label{omegadef}
\omega(\varepsilon)=c_0\cdot\sqrt[3]{\varepsilon}|\ln \varepsilon|^{\frac43}
\end{equation}
 be the error estimate in Corollary~\ref{PLBstab}.

The main goal of this section is to prove

\begin{lemma}
If $0<\varepsilon<\frac1{c\,n^4}$ for suitable absolute constant $c>1$, then
\begin{equation}
\label{levelset-difference}
\begin{array}{rcl}
\int_0^\infty\left|\,|\Phi_t|-|\Omega_t|\,\right|\,dt&\leq & 97\sqrt{n}\, \sqrt{\omega(\varepsilon)} \\[1ex]
\int_0^\infty\left|\,|\Psi_t|-|\Omega_t|\,\right|\,dt&\leq & 97\sqrt{n}\, \sqrt{\omega(\varepsilon)}.
\end{array}
\end{equation}
\end{lemma}
\proof 
The absolute constant $c>1$ in the condition $0<\varepsilon<\frac1{c\,n^4}$ is defined in a way such that ({\it cf} \eqref{omegadef})
\begin{equation}
\label{levelset-differencecdef0}
\omega(\varepsilon)<\frac1{4\cdot 24^2n},
\end{equation}
which inequality is equivalent with \eqref{levelset-differencecdef} below.
 
We observe that $\Phi_t,\Psi_t,\Omega_t$
are  convex bodies, and
$F(t),G(t),H(t)$ are decreasing and log-concave,
and $F,G$ are probability distributions
on $[0,\infty)$ by \eqref{FGint1}. Since $\int_0^\infty H=\int_{\R^n}  h\leq 1+\varepsilon$,
it follows from Corollary~\ref{PLBstab} that
there exists some $b>0$ such that
\begin{equation}
\label{FHGH}
\begin{array}{rcl}
\int_0^\infty|bF(bt)-H(t)|\,dt&\leq&\omega(\varepsilon)\\
\int_0^\infty|b^{-1}G(b^{-1}t)-H(t)|\,dt&\leq &\omega(\varepsilon).
\end{array}
\end{equation}
We may assume that $b\geq 1$.

For $t>0$, let
\begin{eqnarray*}
\widetilde{\Phi}_t&=&
b^{\frac{1}n}\Phi_{bt}
\mbox{ \ \ if $\widetilde{\Phi}_t\neq\emptyset$} \\
\widetilde{\Psi}_t&=&b^{\frac{-1}n}\Psi_{b^{-1}t} 
\mbox{ \ \ if $\widetilde{\Psi}_t\neq\emptyset$}.
\end{eqnarray*}
These sets satisfy $|\widetilde{\Phi}_t|=bF(bt)$, $|\widetilde{\Psi}_t|=b^{-1}G(b^{-1}t)$ and
\begin{eqnarray}
\label{tildePhiH}
\int_0^\infty|\,|\widetilde{\Phi}_t|-H(t)|\,dt&\leq &\omega(\varepsilon)\\
\label{tildePsiH}
\int_0^\infty|\,|\widetilde{\Psi}_t|-H(t)|\,dt&\leq &
\omega(\varepsilon).
\end{eqnarray}
In addition, \eqref{sections} yields that
if $\widetilde{\Phi}_t\neq\emptyset$ and $\widetilde{\Psi}_t\neq\emptyset$, then
\begin{equation}
\label{widetildeOmega}
 \mbox{$\frac12$}(b^{\frac{-1}n}\widetilde{\Phi}_t
+b^{\frac{1}n}\widetilde{\Psi}_t)\subset\Omega_t.
\end{equation}

We dissect $[0,\infty)$ into $I$ and $J$,
where $t\in I$, if $\frac34\,H(t)<|\widetilde{\Phi}_t|<\frac54\,H(t)$
and $\frac34\,H(t)<|\widetilde{\Psi}_t|<\frac54\,H(t)$, and
$t\in J$ otherwise. 
For $J$,
since $\varepsilon<\frac1{c\,n^4}$ and we choose $c>1$ in a way such that \eqref{levelset-differencecdef}
holds, 
\eqref{tildePhiH} and \eqref{tildePsiH} yield that 
\begin{equation}
\label{Jsize}
\int_J H(t)\,dt\leq 4\int_J\left(|\,|\widetilde{\Phi}_t|-H(t)|+
|\,|\widetilde{\Psi}_t|-H(t)|\right)\,dt
\leq 8\omega(\varepsilon)<\frac12.
\end{equation}

Turning to $I$, it follows from the Pr\'ekopa-Leindler inequality
and (\ref{Jsize}) that
\begin{equation}
\label{Isize}
\int_I H(t)\,dt\geq 1-\int_J H(t)\,dt>\frac12.
\end{equation}
For $t\in I$, we define
 $\alpha(t)=|\widetilde{\Phi}_t|/H(t)$ and
$\beta(t)=|\widetilde{\Psi}_t|/H(t)$, and hence
$\frac34<\alpha(t),\beta(t)<\frac54$, and
\eqref{tildePhiH} and \eqref{tildePsiH} imply
\begin{equation}
\label{alphabeta}
\int_0^\infty H(t)\cdot\left(|\alpha(t)-1|+|\beta(t)-1|\right)\,dt\leq 2\omega(\varepsilon).
\end{equation}

In addition, let
\begin{eqnarray*}
\sigma(t)&=&\sigma\left(b^{\frac{-1}n}\widetilde{\Phi}_t,
b^{\frac{1}n}\widetilde{\Psi}_t\right)=
\max\left\{\frac{b^2\beta(t)}{\alpha(t)},
\frac{\alpha(t)}{b^2\beta(t)}\right\}\\
\eta(t)&=&
\frac{(\sigma(t)-1)^2}{32n\sigma(t)^2}
+\frac{n\gamma^*}{\sigma(t)^{\frac1n}}\cdot
A(\widetilde{\Phi}_t,\widetilde{\Psi}_t)^2,
\end{eqnarray*}
where $\gamma^*$ comes from Theorem~\ref{Maggi} and \eqref{Maggiprod}.
For $t\in I$, we observe that $\sqrt{\alpha(t)\cdot \beta(t)}\geq 1-\max\{0,1-\alpha(t)\}-\max\{0,1-\beta(t)\}>\frac12$
holds by $\alpha(t),\beta(t)>\frac34$, and hence
(\ref{Maggiprod}) and (\ref{widetildeOmega}) yields
that
\begin{eqnarray*}
H(t)&\geq& 
\sqrt{\left|b^{\frac{-1}n}\widetilde{\Phi}_t\right|\cdot \left|b^{\frac{1}n}\widetilde{\Psi}_t\right|}(1+\eta(t))
=H(t)\cdot\sqrt{\alpha(t)\cdot \beta(t)}(1+\eta(t))\\
&\geq& H(t)\cdot\left(1-\max\{0,1-\alpha(t)\}-\max\{0,1-\beta(t)\}\right)
(1+\eta(t))\\
&\geq& H(t)\cdot(1-|\alpha(t)-1|-|\beta(t)-1|+\mbox{$\frac12$}\,\eta(t)).
\end{eqnarray*}
In particular, (\ref{alphabeta}) implies
\begin{equation}
\label{Ieta}
\int_I  H(t)\cdot\eta(t)\,dt\leq 4\omega(\varepsilon).
\end{equation}

Next we estimate $b\geq 1$ from above (see \eqref{best}). 
First we claim that if $t\in I$, then
\begin{equation}
\label{besttI}
|\alpha(t)-1|+|\beta(t)-1|+\eta(t)\geq\frac{(b-1)^2}{32n b^2}.
\end{equation}
If $\alpha(t)< b\beta(t)$, then $\sigma(t)> b$, and hence
$$
\eta(t)>\frac{(b-1)^2}{32n b^2}.
$$
Thus let $\alpha(t)\geq b\beta(t)$ in \eqref{besttI}. 
If $\beta(t)\geq 1$, then
$$
|\alpha(t)-1|+|\beta(t)-1|\geq |\alpha(t)-1|\geq b-1\geq
\frac{b-1}{b},
$$
if $1/b\leq \beta(t)\leq 1$, then
$$
|\alpha(t)-1|+|\beta(t)-1|\geq b\beta(t)-1+1-\beta(t)\geq 
\frac{b-1}{b},
$$
and if $\beta(t)\leq 1/b$, then
$$
|\alpha(t)-1|+|\beta(t)-1|\geq|\beta(t)-1|\geq 1-\frac1b=
\frac{b-1}{b},
$$
completing the proof of \eqref{besttI}.

We deduce from (\ref{Isize}), (\ref{alphabeta}), (\ref{Ieta}) and \eqref{besttI} that
\begin{eqnarray*}
\frac{(b-1)^2}{64n b^2}&\leq & \int_I H(t)\cdot\frac{(b-1)^2}{32n b^2}\,dt\\
&\leq & \int_0^\infty H(t)\cdot
\left(\eta(t)+|\alpha(t)-1|+|\beta(t)-1|\right)\,dt\leq 6\,\omega(\varepsilon),
\end{eqnarray*}
and hence
$$
\frac{b-1}{ b}\leq 24\sqrt{n}\sqrt{\omega(\varepsilon)}.
$$
We have choosen $c>1$ in the condition $\varepsilon<\frac1{c\,n^4}$ large enough such that
({\it cf.} \eqref{levelset-differencecdef0})
\begin{equation}
\label{levelset-differencecdef}
24\sqrt{n}\sqrt{\omega(\varepsilon)}<\frac12.
\end{equation}

Since $\frac{b-1}{b}>\frac12$ if $b>2$, we conclude that
\begin{equation}
\label{best}
b-1\leq 48\sqrt{n}\,\sqrt{\omega(\varepsilon)}.
\end{equation}

Next we claim that
\begin{equation}
\label{PhitildePhi}
\begin{array}{rcl}
\int_0^\infty\left|\,|\Phi_t|-|\widetilde{\Phi}_t|\,\right|\,dt&\leq & 96\sqrt{n}\, \sqrt{\omega(\varepsilon)} \\[1ex]
\int_0^\infty\left|\,|\Psi_t|-|\widetilde{\Psi}_t|\,\right|\,dt&\leq &  96\sqrt{n}\, \sqrt{\omega(\varepsilon)}.
\end{array}
\end{equation}
Since $|\Phi_{bt}|\leq |\Phi_t|$, we have
\begin{eqnarray*}
\int_0^\infty\left||\Phi_t|-|\widetilde{\Phi}_t|\right|\,dt&=&\int_0^\infty\left||\Phi_t|-b|\Phi_{bt}|\right|\,dt\\
&\leq &\int_0^\infty\left||\Phi_t|-b|\Phi_{t}|\right|\,dt+b \int_0^\infty\left||\Phi_t|-|\Phi_{bt}|\right|\,dt\\
&= & (b-1)+
b\int_0^\infty|\Phi_t|-|\Phi_{bt}|\,dt\\
&=& 2(b-1)\leq 96\sqrt{n}\,\sqrt{\omega(\varepsilon)}.
\end{eqnarray*}
Similarly, $|\Psi_t|\leq |\Psi_{b^{-1}t}|$, and hence
\begin{eqnarray*}
\int_0^\infty\left||\Psi_t|-|\widetilde{\Psi}_t|\right|\,dt&=&\int_0^\infty\left||\Psi_t|-b^{-1}|\Psi_{b^{-1}t}|\right|\,dt\\
&\leq &\int_0^\infty\left||\Psi_t|-b^{-1}|\Psi_{t}|\right|\,dt+
b^{-1} \int_0^\infty\left||\Psi_t|-|\Psi_{b^{-1}t}|\right|\,dt\\
&= & (1-b^{-1})+
b^{-1}\int_0^\infty|\Psi_{b^{-1}t}|-|\Psi_t|\,dt\\
&=& 2(1-b^{-1})\leq 96\sqrt{n}\,\sqrt{\omega(\varepsilon)},
\end{eqnarray*}
proving (\ref{PhitildePhi}).

We conclude \eqref{levelset-difference} from combining \eqref{tildePhiH}, \eqref{tildePsiH} and  (\ref{PhitildePhi}).
\proofbox

As a first consequence of \eqref{levelset-difference}, we verify the following.

\begin{coro}
\label{fogo212}
There exists an absolute constant $c>1$ such that
if $0<\varepsilon<(cn)^{-n}$, then
$\frac12<M_f/M_g<2$ and $\frac12<M_f/M_h<2$.
\end{coro}
\proof First we prove the estimate about $M_g$ using that
\eqref{levelset-difference} yields
\begin{equation}
\label{levelset-differencePhiPsi}
\int_0^\infty\left|\,|\Phi_t|-|\Psi_t|\,\right|\,dt\leq 194 \sqrt{\omega(\varepsilon)}.
\end{equation}

We may assume that $1=M_f\geq M_g$. Since $|\Psi_t|=0$ if $t>M_g$, we deduce first from 
\eqref{levelset-differencePhiPsi}, and then  from \eqref{closetomax} 
and $k!<(\frac{k}{e})^k\sqrt{2\pi (k+1)}$ that
\begin{eqnarray*}
194\sqrt{\omega(\varepsilon)}&\geq& 
\int_{M_g}^1|\Phi_t|\,dt\geq \frac1{2\cdot n!} \int_{M_g}^1 (1-t)^n\,dt\\
&=&\frac1{2\cdot n!}\frac{(1-M_g)^{n+1}}{n+1}>
\frac{e^{n+1}}{2(n+1)^{n+1}\sqrt{2\pi(n+2)}}\cdot (1-M_g)^{n+1},
\end{eqnarray*}
and hence 
$$
1-M_g< c_1n\omega(\varepsilon)^{\frac1{2(n+1)}}
$$
for an absolute constant $c_1>0$. In particular, we deduce from
\eqref{omegadef} that for some absolute constant $c>1$,
if $0<\varepsilon<(cn)^{-n}$, then $M_g>\frac12$.

The proof of $\frac12<M_f/M_h<2$ is analogous based directly on \eqref{levelset-difference}.
\proofbox

\section{Proof of Theorem~\ref{PLhstab12}}
\label{sechalf}

We use the notation set up in Section~\ref{seclevelsets} with the additional assumption $f(o)=1$, and hence
\begin{equation}
\label{fo1go}
f(o)=M_f=1\mbox{ and }g(o)=M_g.
\end{equation}
First we assume that
\begin{equation}
\label{epscondTh31}
\varepsilon<c^{-n}n^{-n}
\end{equation}
for suitably large absolute constant $c>1$.
According to \eqref{fo1go}, \eqref{epscondTh31} and Corollary~\ref{fogo212}, we have
\begin{equation}
\label{g0bound}
\begin{array}{rcl}
\frac12<&g(o)=M_g&<2\\[0.5ex]
\frac12<&M_h&<2.
\end{array}
\end{equation}

We assume that 
$\R^n$ is a linear subspace of $\R^{n+1}$, and
 write $u_0$ to denote the $(n+1)$th basis vector in $\R^{n+1}$ orthogonal to $\R^n$.
We set 
\begin{equation}
\label{xidef}
\xi=\frac{\sqrt[6]{\omega(\varepsilon)}}{|\ln\omega(\varepsilon)|^{\frac12}}
\end{equation}
where \eqref{epscondTh31} ensures that
\begin{equation}
\label{xiest}
\xi<\frac{e^{-4(n-1)}}2
\mbox{ \ and \ }6e\xi \cdot|\ln \xi|^n<\frac12.
\end{equation}
 In particular,
using the substitution $s=\ln t$, it follows from $M_f=1$ (see \eqref{fo1go}), 
$\frac12<M_g,M_h<2$ (see \eqref{g0bound}) and \eqref{philessthans}) and \eqref{xiest} that
\begin{eqnarray}
\label{intxiPhitlower}
\int_{\xi}^{1}|\Phi_t|\,dt=\int_{\xi}^{M_f}|\Phi_t|\,dt&>&1-2e\cdot \xi\cdot|\ln \xi|^n>\frac12\\
\label{intxiPsitlower}
\int_{\xi}^{2}|\Psi_t|\,dt=\int_{\xi}^{M_g}|\Psi_t|\,dt&>&
1- 3e\cdot \frac{\xi}{M_g}\cdot \left|\ln \frac{\xi}{M_g}\right|^n\\
\nonumber
&>&1-6e\xi\cdot|\ln \xi|^n>\frac12\\
\label{intxiOmegatlower}
\int_{\xi}^{2}|\Omega_t|\,dt=\int_{\xi}^{M_h}|\Omega_t|\,dt&>&
1-6e\xi\cdot|\ln \xi|^n>\frac12.
\end{eqnarray}

We consider the following convex bodies
in $\R^{n+1}$:
\begin{eqnarray}
\label{Kdef}
K=K_{\xi,f}&=&\{x+u_0\ln t:\,x\in\Phi_\xi\mbox{ and }\xi\leq t\leq f(x)\}\\
\label{Cdef}
C=C_{\xi,g}&=&\{x+u_0\ln t:\,x\in\Psi_\xi\mbox{ and }\xi\leq t\leq g(x)\}\\
\label{Ldef}
L=L_{\xi,h}&=&\{x+u_0\ln t:\,x\in\Omega_\xi\mbox{ and }\xi\leq t\leq h(x)\}.
\end{eqnarray}

We write $V(\cdot)$ to denote volume ($(n+1)$-dimensional Lebesgue measure) in $\R^{n+1}$. 
It follows from \eqref{intxiPhitlower} and  \eqref{intxiPsitlower} that
\begin{equation}
\label{VKClow}
\begin{array}{rcl}
V(K)&=&
\int_{\ln \xi}^{0}|\Phi_{e^s}|\,ds
=\int_{\xi}^{1}|\Phi_t|\cdot\frac1t\,dt\geq
\int_{\xi}^{1}|\Phi_t|\,dt>\frac12,\\
V(C)&=&
\int_{\ln \xi}^{\ln 2}|\Psi_{e^s}|\,ds=\int_{\xi}^{2}|\Psi_t|\cdot\frac1t\,dt\geq
\int_{\xi}^{2}|\Psi_t|\cdot\frac12\,dt>\frac14.
\end{array}
\end{equation}
Since $K$ is contained in a right cylinder whose base is a translate of 
$\Phi_{\xi}$ and height is $|\ln \xi|$,
and $C$ is contained in a right cylinder whose base is a copy of 
$\Psi_{\xi}$ and height is $|\ln \xi|+\ln 2<2|\ln \xi|$, we deduce from 
\eqref{philessthans0} that
\begin{equation}
\label{VKCupp}
\begin{array}{rcl}
V(K)&\leq&
\frac{2}{n!} \cdot |\ln \xi|^{n+1},\\[0.5ex]
V(C)&\leq&
\frac{4}{n!} \cdot |\ln \xi|^{n+1}.
\end{array}
\end{equation}

It follows from (\ref{levelset-difference}), $f(o)=1$, (\ref{g0bound}) and using the substitution $s=\ln t$  that
\begin{eqnarray}
\nonumber
|V(K)-V(L)|&=&
\left|\int_{\ln \xi}^{\ln 2}(|\Phi_{e^s}|-|\Omega_{e^s}|)\,ds\right|
=\left|\int_{\xi}^{2}(|\Phi_t|-|\Omega_t|)\cdot\frac1t\,dt\right|\\
\label{VKVL}
&\leq &\frac1{\xi}\int_\xi^{2}\left| |\Phi_t|-|\Omega_t|\right|\,dt
\leq 97\sqrt{n}\,\frac{\sqrt{\omega(\varepsilon)}}{\xi},
\end{eqnarray}
and similarly
\begin{equation}
\label{VCVL}
|V(C)-V(L)|\leq 97\sqrt{n}\,
\frac{\sqrt{\omega(\varepsilon)}}{\xi}.
\end{equation}
Combining \eqref{VKVL} and  \eqref{VCVL} leads to
\begin{equation}
\label{VCVK}
|V(C)-V(K)|\leq 194\sqrt{n}\,
\frac{\sqrt{\omega(\varepsilon)}}{\xi}.
\end{equation}

The condition $h(\frac{x+y}2)\geq \sqrt{f(x)g(y)}$ for
$x,y\in\mathbb{R}^n$ in Theorem~\ref{PLhstab12} implies that
\begin{equation}
\label{K+CL}
\mbox{$\frac12$}\,K+\mbox{$\frac12$}\,C\subset L.
\end{equation}

\begin{lemma}
\label{KCLtranstales}
Assuming the condition $\varepsilon<c^{-n}n^{-n}$ as in \eqref{epscondTh31},
there exist $w\in\R^n$ and absolute constant $\gamma>1$ such that
\begin{eqnarray*}
V(K\Delta (L-w))&\leq& \gamma n^8\cdot\frac{\sqrt[4]{\omega(\varepsilon)}}{\sqrt{\xi}}\cdot |\ln \xi|^{\frac{n+1}2}\\
V(C\Delta (L+w))&\leq&\gamma n^8\cdot\frac{\sqrt[4]{\omega(\varepsilon)}}{\sqrt{\xi}}\cdot |\ln \xi|^{\frac{n+1}2}.
\end{eqnarray*}
\end{lemma}
\proof First we verify a slightly weaker statement; namely, we allow to choose the translation
vectors from $\R^{n+1}$, not only from $\R^n$. More precisely, we claim that there exist
$\tilde{w}\in\R^{n+1}$ and absolute constant $\gamma>1$ such that
\begin{eqnarray}
\label{KLtranslatetildez}
V(K\Delta (L-\tilde{w}))&\leq&\frac{\gamma n^8}3\cdot\frac{\sqrt[4]{\omega(\varepsilon)}}{\sqrt{\xi}}\cdot |\ln \xi|^{\frac{n+1}2}\\
\label{CLtranslatetildew}
V(C\Delta (L+\tilde{w}))&\leq&\frac{\gamma n^8}3\cdot
\frac{\sqrt[4]{\omega(\varepsilon)}}{\sqrt{\xi}}\cdot |\ln \xi|^{\frac{n+1}2}.
\end{eqnarray}
To prove  \eqref{KLtranslatetildez} and \eqref{CLtranslatetildew}, let
 us consider a homothetic copy $K_0\subset K$ of $K$ and homothetic copy $C_0\subset C$ of $C$
with $V(K_0)=V(C_0)=\min\{V(K),V(C)\}$, and hence either $K=K_0$ or $C=C_0$, and
\begin{equation}
\label{K0+C0L}
\mbox{$\frac12$}\,K_0+\mbox{$\frac12$}\,C_0\subset L.
\end{equation}
It also follows from \eqref{VKVL} and  \eqref{VCVL}, and from the fact that
 either $K=K_0$ or $C=C_0$,  that
\begin{equation}
\label{VLleVK0}
V(L)-V(K_0)\leq 97\sqrt{n}\cdot \frac{\sqrt{\omega(\varepsilon)}}{\xi\,V(K_0)}\cdot V(K_0) 
\end{equation}
where \eqref{epscondTh31} (provided $\tilde{c}$ is large enough), \eqref{VKClow} and 
$\xi=\sqrt[6]{\omega(\varepsilon)}/|\ln\omega(\varepsilon)|^{\frac12}$
yield that 
$\frac{\sqrt{\omega(\varepsilon)}}{\xi\,V(K_0)}<2\omega(\varepsilon)^{\frac13}|\ln\omega(\varepsilon)|^{\frac12}$
is small enough to apply Lemma~\ref{VKC1L} to $K_0$, $C_0$ and $L$. In particular, we deduce from \eqref{K0+C0L}, \eqref{VLleVK0} and Lemma~\ref{VKC1L}  that there exist
$\tilde{w}\in\R^{n+1}$ and an absolute constant $\gamma_0>1$ such that
\begin{eqnarray*}
V(K_0\Delta (L-\tilde{w}))&\leq& \gamma_0n^8\cdot \frac{\sqrt[4]{\omega(\varepsilon)}}{\sqrt{\xi\,V(K_0)}}\cdot V(K_0)=
\gamma_0n^8\frac{\sqrt[4]{\omega(\varepsilon)}}{\sqrt{\xi}}\cdot \sqrt{V(K_0)}\\
V(C_0\Delta (L+\tilde{w}))&\leq&\gamma_0n^8\cdot \frac{\sqrt[4]{\omega(\varepsilon)}}{\sqrt{\xi}}\cdot \sqrt{V(K_0)}.
\end{eqnarray*}
In turn, it follows from \eqref{VKClow}, \eqref{VCVK} and the properties of $K_0$ and $C_0$  that
\begin{eqnarray*}
V(K\Delta (L-\tilde{w}))&\leq&(\gamma_0+1)n^8 \cdot 
\frac{\sqrt[4]{\omega(\varepsilon)}}{\sqrt{\xi}}\cdot \sqrt{V(K)}\\
V(C\Delta (L+\tilde{w}))&\leq&(\gamma_0+1)n^8\cdot \frac{\sqrt[4]{\omega(\varepsilon)}}{\sqrt{\xi}}\cdot \sqrt{V(K)}.
\end{eqnarray*}
Therefore, \eqref{VKCupp} implies \eqref{KLtranslatetildez} and \eqref{CLtranslatetildew}.

Next we verify that if $w\in\R^n$
satisfies that $\tilde{w}=w+pu_0$ for $p\in \R$, then, then we have
\begin{eqnarray}
\label{Kztildezvol}
V(K\Delta (L-w))&\leq & 3V(K \Delta (L-\tilde{w}))\\
\label{Cwtildezvol}
V(C\Delta (L+w))&\leq & 3V(C \Delta (L+\tilde{w})).
\end{eqnarray}  
For \eqref{Kztildezvol}, we may assume that $p\neq 0$, and we distinguish two cases.

If $p<0$, then let
$$
K_{(p)}=\{x\in K:\,\ln\xi\leq \langle x,u_0\rangle<|p|+\ln\xi\}\subset K\backslash (L-\tilde{w}).
$$
Using the fact that $\Phi_t$ is decreasing as a set as $t>0$ increases and the Fubini theorem, we have
$$
V(K\Delta (K+|p|u_0))=2V(K_{(p)})\leq 2V(K\Delta (L-\tilde{w}));
$$
therefore, the triangle inequality for the symmetric difference metric implies
\begin{eqnarray*}
V(K\Delta (L-w))&=& V((K+|p|\,u_0)\Delta (L-\tilde{w}))\\
&\leq& V((K+|p|\,u_0)\Delta K)+V(K \Delta (L-\tilde{w}))\\
&\leq & 3V(K \Delta (L-\tilde{w})).
\end{eqnarray*}

Similarly, if $p>0$, then we consider
$$
L_{(p)}=\{x\in L:\,\ln\xi \leq \langle x,u_0\rangle<p+\ln\xi\}
$$
satisfying
$$
L_{(p)}+\tilde{w}\subset (L+\tilde{w})\backslash K.
$$
Using that $\Omega_t$ is decreasing as a set as $t>0$ increases, we deduce
$$
V((L+\tilde{w})\Delta (L+w))=2V(L_{(p)})\leq 2V(K\Delta (L+\tilde{w}));
$$
therefore, the triangle inequality for the symmetric difference metric implies
$$
V(K\Delta (L+w))\leq V(K \Delta (L+\tilde{w}))+V((L+\tilde{w})\Delta (L+w) )
\leq  3V(K \Delta (L+\tilde{w})),
$$
completing the proof of \eqref{Kztildezvol}.
The argument for \eqref{Cwtildezvol} is similar.

 In turn, combining \eqref{KLtranslatetildez}, \eqref{CLtranslatetildew},
\eqref{Kztildezvol} and \eqref{Cwtildezvol} yields Lemma~\ref{KCLtranstales}.
\proofbox

\noindent{\bf Proof of Theorem~\ref{PLhstab12}}
We may assume that $f$ and $g$ are log-concave probability distributions satisfying
$$
f(o)=M_f=1\mbox{ \ and \ } g(o)=M_g.
$$
In particular, \eqref{g0bound} implies
$\frac12<M_g,M_h<2$.
We consider the convex bodies $K,C,L\subset \R^{n+1}$ defined by \eqref{Kdef}, \eqref{Cdef}
and \eqref{Ldef}, and let $w\in\R^n$ be provided by Lemma~\ref{KCLtranstales}.

For \eqref{philessthans}, we have $\frac12<M_f,M_h<2$, and hence $\frac{\xi}{M_f},\frac{\xi}{M_g}<2\xi$
and both $1+\frac1{M_f}$ and $1+\frac1{M_h}$ are at most $3$.
For the functions $f,g,h$, 
it follows from $n\geq 2$, \eqref{philessthans} (compare the condition \eqref{xiest}),  Lemma~\ref{KCLtranstales},
$\xi=\sqrt[6]{\omega(\varepsilon)}/|\ln\omega(\varepsilon)|^{\frac12}$ and the $\gamma$ of 
Lemma~\ref{KCLtranstales} that
\begin{eqnarray*}
\int_{\R^n}|f(x)-h(x-w)|\,dx&=& \int_0^2 \left|\Phi_t\Delta (\Omega_t-w) \right|\,dt\\
&\leq &\int_{\xi}^2 \left|\Phi_t\Delta (\Omega_t-w) \right|\,dt+
\int_0^{\xi} \left|\Phi_t\right|\,dt+
\int_0^{\xi} \left|\Omega_t\right|\,dt\\
&\leq &2\int_{\xi}^2 \left|\Phi_t\Delta (\Omega_t-w) \right|\cdot \frac1t\,dt+
\int_0^{\xi} \left|\Phi_t\right|\,dt+
\int_0^{\xi} \left|\Omega_t\right|\,dt\\
&= &2V(K\Delta(L-w))+
\int_0^{\xi} \left|\Phi_t\right|\,dt+
\int_0^{\xi} \left|\Omega_t\right|\,dt\\
&\leq &\gamma n^8\cdot\frac{\sqrt[4]{\omega(\varepsilon)}}{\sqrt{\xi}}\cdot |\ln \xi|^{\frac{n+1}2}+
2\cdot 3\cdot (2\xi)\cdot |\ln (2\xi)|^n \\
&\leq & 2\gamma n^8\cdot\sqrt[6]{\omega(\varepsilon)}\cdot |\ln \omega(\varepsilon)|^{n-\frac14}.
\end{eqnarray*}

Similarly, we have
$$
\int_{\R^n}|g(x)-h(x+w)|\,dx\leq 2\gamma n^8\cdot \sqrt[6]{\omega(\varepsilon)}
\cdot |\ln \omega(\varepsilon)|^{n-\frac14}.
$$
Since $\omega(\varepsilon)=c_0\sqrt[3]{\varepsilon}|\ln \varepsilon|^{\frac43}$ 
for an absolute constan $c_0>1$ ({\it cf.} \eqref{omegadef}), we conclude that
\begin{eqnarray*}
\int_{\R^n}|f(x)-h(x-w)|\,dx&\leq& \gamma_0 n^8\cdot \sqrt[18]{\varepsilon}\cdot|\log\varepsilon|^{n} \\
\int_{\R^n}|g(x)-h(x+w)|\,dx&\leq&\gamma_0 n^8\cdot \sqrt[18]{\varepsilon}\cdot|\log\varepsilon|^{n} 
\end{eqnarray*}
for an absolute constant $\gamma_0>1$, proving Theorem~\ref{PLhstab12}. 
\proofbox

\section{A version of Theorem~\ref{PLhstab0} when $\varepsilon$ is small}
\label{seclambda}

The goal of this section is to prove the following  version of Theorem~\ref{PLhstab0}.

\begin{theo}
\label{PLhstab}
For some absolute constant $c>1$, if $\tau\in(0,\frac12]$, $\tau\leq \lambda\leq 1-\tau$, 
$h,f,g:\,\R^n\to [0,\infty)$ are integrable  such that
   $h((1-\lambda)x+\lambda\,y)\geq f(x)^{1-\lambda}g(y)^\lambda$ for
$x,y\in\mathbb{R}^n$, $h$ is log-concave and
$$
\int_{\R^n}  h\leq (1+\varepsilon) \left(\int_{\R^n}f\right)^{1-\lambda}\left(\int_{\R^n}g\right)^\lambda
$$
for $\varepsilon\in(0,\tau \varepsilon_0)$ for $\varepsilon_0=c^{-n}n^{-n}$, then
there exists $w\in\R^n$ such that setting $a=\int_{\R^n}f/\int_{\R^n}g$, we have
\begin{eqnarray*}
\int_{\R^n}|f(x)-a^\lambda h(x-\lambda\,w)|\,dx&\leq &cn^8\sqrt[18]{\frac{\varepsilon}{\tau}}\cdot
\left|\log\frac{\varepsilon}{\tau}\right|^{n}  
\int_{\R^n}f \\
\int_{\R^n}|g(x)-a^{-(1-\lambda)}h(x+(1-\lambda)w)|\,dx&\leq &cn^8\sqrt[18]{\frac{\varepsilon}{\tau}}\cdot
\left|\log\frac{\varepsilon}{\tau}\right|^{n}    
\int_{\R^n}g.
\end{eqnarray*}
\end{theo}

For a bounded measurable function $f:\,\R^n\to[0,\infty)$, the log-concave hull $\tilde{f}:\,\R^n\to[0,\infty)$ of is
$$
\tilde{f}(z)=\sup_{z=\sum_{i=1}^k \alpha_i x_i \atop \sum_{i=1}^k \alpha_i=1,\;\forall \alpha_i\geq 0}
\prod_{i=1}^k f(x_i)^{\alpha_i}.
$$
To show that $\tilde{f}$ is log-concave, it is equivalent to prove that if
$\varepsilon,\alpha,\beta\in(0,1)$ and $x,y\in\R^n$, then
\begin{equation}
\label{tildefepslogconcave}
\tilde{f}(\alpha x+\beta y)\geq (1-\varepsilon)\tilde{f}(x)^\alpha\tilde{f}(y)^\beta.
\end{equation}
We observe that there exist $x_1,\ldots,x_k,y_1,\ldots,y_m\in\R^n$ and
$\alpha_1,\ldots,\alpha_k,\beta_1,\ldots,\beta_m\geq 0$ with
$\sum_{i=1}^k \alpha_i=1$, $\sum_{j=1}^m \beta_j=1$,
$x=\sum_{i=1}^k \alpha_i x_i$ and $y=\sum_{i=1}^k \beta_j x_j$ such that
$$
(1-\varepsilon)\tilde{f}(x)\geq \prod_{i=1}^k f(x_i)^{\alpha_i}
\mbox{ \ and \ }
(1-\varepsilon)\tilde{f}(y)\geq \prod_{j=1}^m f(y_j)^{\beta_j}.
$$
Since
\begin{equation}
\label{alphabetaij}
\alpha x+\beta y=\sum_{i=1}^k\sum_{j=1}^m (\alpha_i\beta_j \alpha x_i
+\alpha_i\beta_j \beta y_j) \mbox{ \ where \ }
\sum_{i=1}^k\sum_{j=1}^m (\alpha_i\beta_j \alpha 
+\alpha_i\beta_j \beta)=1,
\end{equation}
we deduce that
\begin{eqnarray*}
\tilde{f}(\alpha x+\beta y)&\geq&
\prod_{i=1}^k\prod_{j=1}^m f(x_i)^{\alpha_i\beta_j \alpha}
f(y_j)^{\alpha_i\beta_j \beta}\\
&=&\left(\prod_{i=1}^k f(x_i)^{\alpha_i}\right)^\alpha
\left(\prod_{j=1}^m f(y_j)^{\beta_j }\right)^{\beta}\\
&\geq & (1-\varepsilon)^\alpha\tilde{f}(x)^\alpha(1-\varepsilon)^\beta\tilde{f}(y)^\beta
=(1-\varepsilon)\tilde{f}(x)^\alpha\tilde{f}(y)^\beta,
\end{eqnarray*}
proving that $\tilde{f}$ is log-concave via \eqref{tildefepslogconcave}.

We note that if $a_0>0$ and $z_0\in\R^n$ and $f_0(z)=a_0f(z-z_0)$, then
\begin{equation}
\label{tildefa0z0}
\widetilde{f_0}(z)=a_0\tilde{f}(z-z_0).
\end{equation}

We prepare the proof of Theorem~\ref{PLhstab}
by the three technical statements
Lemma~\ref{log-concave-below},
Lemma~\ref{log-concave-geometric-mean} and
Lemma~\ref{log-concave-error} about log-concave functions.

\begin{lemma}
\label{log-concave-below}
If  $\lambda\in(0,1)$, $h$ is a log-concave function on $\R^n$ with positive integral, and $f,g:\,\R^n\to[0,\infty)$
are measurable satisfying $\int_{\R^n} f>0$, $\int_{\R^n} g>0$ and  
$h((1-\lambda)x+\lambda y)\geq f(x)^{1-\lambda}g(y)^\lambda$ for
$x,y\in\R^n$, then $f$ and $g$ are bounded, and their log-concave hulls $\tilde{f}$ and $\tilde{g}$ satisfy that
$h((1-\lambda)x+\lambda y)\geq \tilde{f}(x)^{1-\lambda}\tilde{g}(y)^\lambda$ for
$x,y\in\R^n$.
\end{lemma}
{\bf Remark } The Pr\'ekopa-Leindler inequality yields
$\int_{\R^n} \tilde{f}<\infty$ and $\int_{\R^n} \tilde{g}<\infty$.\\
\proof To show that $f$ is bounded, we choose $y_0\in \R^n$ with $g(y_0)>0$. For any $x\in\R^n$, we have
$h((1-\lambda)x+\lambda y_0)\geq f(x)^{1-\lambda}g(y_0)^\lambda$; therefore,
$$
f(x)\leq \frac{h((1-\lambda)x+\lambda y_0)^{\frac{1}{1-\lambda}}}{g(y_0)^{\frac{\lambda}{1-\lambda}}}
\leq \frac{M_h^{\frac{1}{1-\lambda}}}{g(y_0)^{\frac{\lambda}{1-\lambda}}}.
$$
Similar argument yields that $g$ is bounded.

For $x,y\in\R^n$,  it is sufficient to prove that for any $\varepsilon\in(0,1)$, we have
\begin{equation}
\label{htildeftildeg}
h((1-\lambda)x+\lambda y)\geq (1-\varepsilon)\tilde{f}(x)^{1-\lambda}\tilde{g}(y)^\lambda.
\end{equation}
We choose $x_1,\ldots,x_k,y_1,\ldots,y_m\in\R^n$ and
$\alpha_1,\ldots,\alpha_k,\beta_1,\ldots,\beta_m\geq 0$ with
$\sum_{i=1}^k \alpha_i=1$, $\sum_{j=1}^m \beta_j=1$,
$x=\sum_{i=1}^k \alpha_i x_i$ and $y=\sum_{i=1}^k \beta_j x_j$ such that
$$
(1-\varepsilon)\tilde{f}(x)\geq \prod_{i=1}^k f(x_i)^{\alpha_i}
\mbox{ \ and \ }
(1-\varepsilon)\tilde{f}(y)\geq \prod_{j=1}^m f(y_j)^{\beta_j}.
$$
It follows from \eqref{alphabetaij} and the log-concavity of $h$ that
\begin{eqnarray*}
h((1-\lambda)x+\lambda y)&=& h\left(
\sum_{i=1}^k\sum_{j=1}^m \alpha_i\beta_j( (1-\lambda) x_i+ \lambda y_j) \right)\\
&\geq &
\prod_{i=1}^k\prod_{j=1}^m h((1-\lambda) x_i+ \lambda y_j)^{\alpha_i\beta_j }
\geq 
\prod_{i=1}^k\prod_{j=1}^m f(x_i)^{(1-\lambda)\alpha_i\beta_j}g(y_j)^{\lambda \alpha_i\beta_j}\\
&=&\left(\prod_{i=1}^k f(x_i)^{\alpha_i}\right)^{1-\lambda}
\left(\prod_{j=1}^m f(y_j)^{\beta_j }\right)^{\lambda}\\
&\geq & (1-\varepsilon)^{1-\lambda}\tilde{f}(x)^{1-\lambda}(1-\varepsilon)^{\lambda}\tilde{g}(y)^{\lambda}
=(1-\varepsilon)\tilde{f}(x)^{1-\lambda}\tilde{g}(y)^{\lambda},
\end{eqnarray*}
proving \eqref{htildeftildeg}.
\proofbox

\begin{lemma}
\label{log-concave-geometric-mean}
Let $f,g:\,\R^n\to[0,\infty)$ be log-concave and have positive integral.
\begin{description}
\item{(i)} For $\lambda\in[0,1]$, the function $h_\lambda:\,\R^n\to[0,\infty)$ defined by
$$
h_\lambda(z)=\sup_{z=(1-\lambda)x+\lambda y}f(x)^{1-\lambda}g(y)^\lambda
$$
is log-concave, has positive integral, and 
satisfies $h_0=f$ and $h_1=g$.
\item{(ii)} The function $\lambda\mapsto \int_{\R^n}h_\lambda$ is log-concave for $\lambda\in[0,1]$.
\end{description}
\end{lemma}
\proof For (i), readily $h_0=f$ and $h_1=g$. Next let $\lambda\in(0,1)$.
To show the log-concavity of $h_\lambda$, it is sufficient to prove that
 if $z_1,z_2\in\R^n$, $\alpha,\beta>0$ with $\alpha+\beta=1$ and $\varepsilon\in(0,1)$, then
\begin{equation}
\label{hlambdalogconcave}
h_\lambda(\alpha z_1+\beta z_2)\geq (1-\varepsilon)h_\lambda(z_1)^\alpha h_\lambda(z_2)^\beta. 
\end{equation}
We choose $x_1,x_2,y_1,y_2\in\R^n$ such that $z_1=(1-\lambda) x_1+\lambda y_1$, $z_2=
(1-\lambda) x_2+\lambda y_2$
and
$$
f(x_1)^{1-\lambda} g(y_1)^\lambda\geq (1-\varepsilon)h_\lambda(z_1)
\mbox{ \ and \ }
f(x_2)^{1-\lambda} g(y_2)^\lambda\geq (1-\varepsilon)h_\lambda(z_2).
$$
It follows that
$\alpha z_1+\beta z_2=(1-\lambda) (\alpha x_1+\beta x_2)+\lambda (\alpha y_1+\beta y_2)$,
and the log-concavity of $f$ and $g$ yields
\begin{eqnarray*}
h_\lambda(\alpha z_1+\beta z_2)&= &h_\lambda
\big((1-\lambda) (\alpha x_1+\beta x_2)+\lambda (\alpha y_1+\beta y_2)\big)\\
&\geq & f(\alpha x_1+\beta x_2)^{1-\lambda}g(\alpha y_1+\beta y_2)^\lambda\\
&\geq & f(x_1)^{\alpha(1-\lambda)} f(x_2)^{\beta(1-\lambda)}
 g(y_1)^{\alpha(\lambda)} g(y_2)^{\beta(\lambda)}\\
&=& \big(f(x_1)^{1-\lambda}g(y_1)^{\lambda}\big)^{\alpha}
 \big(f(x_2)^{1-\lambda}
  g(y_2)^{\lambda}\big)^{\beta}\\
&\geq & (1-\varepsilon)^{\alpha}h_\lambda(z_1)^{\alpha}
(1-\varepsilon)^{\beta}h_\lambda(z_2)^{\beta}=
(1-\varepsilon)h_\lambda(z_1)^\alpha h_\lambda(z_2)^\beta, 
\end{eqnarray*}
proving \eqref{hlambdalogconcave}, and in turn the log-concavity of $h_\lambda$.

Readily, $\int_{\R^n}h_\lambda>0$. If follows from
Lemma~\ref{log-concave-positive-integral} that $0<M_f,M_g<\infty$, and hence
$$
M=M_{h_\lambda}=M_f^{1-\lambda}M_g^\lambda.
$$
If $t\in (0,M)$ and $h_\lambda(z)>t$, then there exist
$x,y\in \R^n$ such that
\begin{equation}
\label{zxy}
z= (1-\lambda) x+\lambda y
\end{equation}
and $f(x)^{1-\lambda}g(y)^{\lambda}>t$. It follows that
\begin{equation}
\label{zxyfg}
f(x)> \left(\frac{t}{M_g^\lambda}\right)^{\frac1{1-\lambda}}\mbox{ \ and \ }
g(y)> \left(\frac{t}{M_f^{1-\lambda}}\right)^{\frac1{\lambda}}.
\end{equation}
We conclude from \eqref{zxy}, \eqref{zxyfg} and
Lemma~\ref{log-concave-positive-integral} that $h_\lambda(z)>t$ is bounded;
therefore, $h_\lambda$ has positive integral by Lemma~\ref{log-concave-positive-integral}.

Finally, to verify that the function $\lambda\mapsto \int_{\R^n}h_\lambda$ is log-concave for $\lambda\in[0,1]$,
it is enough to prove that if $\lambda_1,\lambda_2\in[0,1]$ and
 $\alpha,\beta>0$ with $\alpha+\beta=1$, then
for $\lambda=\alpha \lambda_1+\beta \lambda_2$, we have
\begin{equation}
\label{inthlambdalogconcave}
\int_{\R^n}h_{\lambda}\geq 
\left(\int_{\R^n}h_{\lambda_1}\right)^\alpha \left(\int_{\R^n}h_{\lambda_2}\right)^\beta. 
\end{equation}
According to the Pr\'ekopa-Leindler inequality Theorem~\ref{PLn}, it is sufficient to show that
if $z=\alpha z_1+\beta z_2$, $z_1,z_2\in\R^n$, then
$$
h_{\lambda}(z)\geq h_{\lambda_1}(z_1)^\alpha
h_{\lambda_2}(z_2)^\beta.
$$
In turn, \eqref{inthlambdalogconcave} is a consequence of the claim that if 
$z=\alpha z_1+\beta z_2$ for $z_1,z_2\in\R^n$ and $\varepsilon\in(0,1)$, then
\begin{equation}
\label{inthlambdalogconcaveeps}
h_{\lambda}(z)\geq (1-\varepsilon) h_{\lambda_1}(z_1)^\alpha
h_{\lambda_2}(z_2)^\beta.
\end{equation}
For $i=1,2$, there exist
$x_i,y_i\in \R^n$ such that
\begin{eqnarray}
\label{zxy12}
z_i&=& (1-\lambda_i) x_i+\lambda_i y_i\\
\label{zxyfg12}
f(x_i)^{1-\lambda_i} g(y_i)^{\lambda_i}&\geq & (1-\varepsilon)h_{\lambda_i}(z_i).
\end{eqnarray}
As $\lambda=\alpha \lambda_1+\beta \lambda_2$ and $z=\alpha z_1+\beta z_2$, we observe that
$1-\lambda=\alpha (1-\lambda_1)+\beta (1-\lambda_2)$ and
\begin{eqnarray*}
z&=&\alpha z_1+\beta z_2=\alpha \big[(1-\lambda_1) x_1+\lambda_1 y_1\big]+
\beta \big[(1-\lambda_2) x_2+\lambda_2 y_2\big]\\
&&(1-\lambda)\cdot \left(\frac{\alpha(1-\lambda_1)}{1-\lambda}\cdot x_1+
\frac{\beta(1-\lambda_2)}{1-\lambda}\cdot x_2\right)+
\lambda\cdot \left(\frac{\alpha \lambda_1}{\lambda}\cdot y_1+
\frac{\beta\lambda_2}{\lambda}\cdot y_2\right).
\end{eqnarray*}
It follows from the log-concavity of $f$ and $g$ and later from \eqref{zxyfg12} that
\begin{eqnarray*}
h_{\lambda}(z)&\geq&
f\left(\frac{\alpha(1-\lambda_1)}{1-\lambda}\cdot x_1+
\frac{\beta(1-\lambda_2)}{1-\lambda}\cdot x_2\right)^{1-\lambda}
g\left(\frac{\alpha \lambda_1}{\lambda}\cdot y_1+
\frac{\beta\lambda_2}{\lambda}\cdot y_2\right)^\lambda\\
&\geq & f(x_1)^{\alpha(1-\lambda_1)} f(x_2)^{\beta(1-\lambda_2)}
g(y_1)^{\alpha\lambda_1} g(y_2)^{\beta\lambda_2}\\
&=&\left(f(x_1)^{1-\lambda_1}g(y_1)^{\lambda_1}\right)^{\alpha}
\left(f(x_2)^{1-\lambda_2}g(y_2)^{\lambda_2}\right)^{\beta}\\
&\geq & (1-\varepsilon)^\alpha h_{\lambda_1}(z_1)^\alpha
(1-\varepsilon)^\beta 
h_{\lambda_2}(z_2)^\beta=(1-\varepsilon) h_{\lambda_1}(z_1)^\alpha
h_{\lambda_2}(z_2)^\beta,
\end{eqnarray*}
proving \eqref{inthlambdalogconcaveeps}, and in turn \eqref{inthlambdalogconcave}.
\proofbox

\begin{lemma}
\label{log-concave-error}
For fixed $\lambda\in(0,1)$, if $\eta\in(0,2\cdot\min\{1-\lambda,\lambda\})$ and $\varphi$ is a log-concave function on $[0,1]$
satisfying  $\varphi(\lambda)\leq (1+\eta)\varphi(0)^{1-\lambda}\varphi(1)^\lambda$, then
$$
\varphi\left(\mbox{$\frac12$}\right)\leq \left(1+\frac{\eta}{\min\{1-\lambda,\lambda\}}\right)\sqrt{\varphi(0)\varphi(1)}
$$
\end{lemma}
\proof We may assume that $0<\lambda<\frac12$, and hence $\lambda=(1-2\lambda)\cdot 0+2\lambda\cdot\frac12$,
$\varphi(\lambda)\leq (1+\eta)\varphi(0)^{1-\lambda}\varphi(1)^\lambda$
and the log-concavity of $\varphi$ yield
$$
(1+\eta)\varphi(0)^{1-\lambda}\varphi(1)^\lambda\geq \varphi(\lambda)\geq
\varphi(0)^{1-2\lambda}\varphi\left(\mbox{$\frac12$}\right)^{2\lambda}.
$$
Thus $(1+\eta)^{\frac1{2\lambda}}\leq e^{\frac{\eta}{2\lambda}}\leq 1+\frac{\eta}{\lambda}$ implies
$$
\varphi\left(\mbox{$\frac12$}\right)\leq (1+\eta)^{\frac1{2\lambda}}\sqrt{\varphi(0)\varphi(1)}\leq
\left(1+\frac{\eta}{\lambda}\right)\sqrt{\varphi(0)\varphi(1)}.
\mbox{ \ \proofbox}
$$

\noindent{\bf Proof of Theorem~\ref{PLhstab}: }
For the $\lambda$ in Theorem~\ref{PLhstab}, we may assume that $0<\lambda\leq \frac12$,
and hence $\min\{1-\lambda,\lambda\}=\lambda$.

For suitable $d,e>0$ and $w\in \R^n$, we may replace $f(z)$ by $d\cdot f(z-w)$, $g(z)$ by
$e\cdot g(z+w)$ and $h(z)$ by $d^{1-\lambda}e^\lambda h(z+(2\lambda-1)w)$
 where $e$ and $d$ will be  defined by \eqref{tildeftildeg1} below, and  $w$ will be defined by 
\eqref{z0fh12} and \eqref{z0gh12}.

Let $\tilde{f}$ and $\tilde{g}$ be the log-concave hulls of $f$ and $g$; therefore,
Lemma~\ref{log-concave-below} yields 
\begin{equation}
\label{hlarger-tildefg}
h((1-\lambda)x+\lambda y)\geq \tilde{f}(x)^{1-\lambda}\tilde{g}(y)^\lambda 
\mbox{ \ for $x,y\in\R^n$}. 
\end{equation}
We may assume by \eqref{tildefa0z0} that
\begin{equation}
\label{tildeftildeg1}
\int_{\R^n}\tilde{f}=\int_{\R^n}\tilde{g}=1.
\end{equation}
It follows from Lemma~\ref{log-concave-geometric-mean} that
$$
h_t(z)=\sup_{z=(1-t)x+t y}\tilde{f}(x)^{1-t}\tilde{g}(y)^t
$$
satisfies that
$$
\varphi(t)=\int_{\R^n}h_t
$$
is log-concave on $[0,1]$. 
In particular, \eqref{tildeftildeg1} implies that
\begin{equation}
\label{phi01}
\varphi(0)=\varphi(1)=1.
\end{equation}

It follows from \eqref{hlarger-tildefg}, \eqref{tildeftildeg1}, the Pr\'ekopa-Leindler inequality Theorem~\ref{PLn} and
the  conditions in Theorem~\ref{PLhstab} that
\begin{eqnarray}
\nonumber
1&=&\left(\int_{\R^n}\tilde{f}\right)^{1-\lambda}\left(\int_{\R^n}\tilde{g}\right)^\lambda\leq 
\int_{\R^n}h_\lambda\leq
\int_{\R^n}h\leq (1+\varepsilon)\left(\int_{\R^n}f\right)^{1-\lambda}\left(\int_{\R^n}g\right)^\lambda \\
\label{PLhstabcor-cond2}
&\leq&
 (1+\varepsilon)\left(\int_{\R^n}\tilde{f}\right)^{1-\lambda}\left(\int_{\R^n}\tilde{g}\right)^\lambda=1+\varepsilon.
\end{eqnarray}
On the other hand,  \eqref{PLhstabcor-cond2}, Lemma~\ref{log-concave-error} 
and finally \eqref{phi01} yield that
$$
\int_{\R^n}h_{1/2}=\varphi\left(\frac12\right)\leq 
\left(1+\frac{\varepsilon}{\lambda}\right)\sqrt{\varphi(0)\varphi(1)}
=1+\frac{\varepsilon}{\lambda}.
$$
In turn, we deduce from Theorem~\ref{PLhstab12} that
there exists $w\in\R^n$ such that
\begin{eqnarray}
\label{z0fh12}
\int_{\R^n}|\tilde{f}(z)-h_{1/2}(z+w)|\,dz&\leq &\tilde{c} n^8\sqrt[18]{\frac{\varepsilon}{\lambda}}
\cdot \left|\log\frac{\varepsilon}{\lambda}\right|^{n},\\
\label{z0gh12}
\int_{\R^n}|\tilde{g}(z)-h_{1/2}(z-w)|\,dz&\leq &\tilde{c} n^8\sqrt[18]{\frac{\varepsilon}{\lambda}}
\cdot \left|\log\frac{\varepsilon}{\lambda}\right|^{n}   .
\end{eqnarray}
Replacing $f(z)$ by $f(z-w)$ and $g(z)$ by
$g(z+w)$ ({\it cf.} \eqref{tildefa0z0}), the function $h_{1/2}$ does not change, and we have
\begin{eqnarray}
\label{fh12}
\int_{\R^n}|\tilde{f}-h_{1/2}|&\leq &\tilde{c} n^8\sqrt[18]{\frac{\varepsilon}{\lambda}}
\cdot \left|\log\frac{\varepsilon}{\lambda}\right|^{n} , \\
\label{gh12}
\int_{\R^n}|\tilde{g}-h_{1/2}|&\leq &\tilde{c} n^8\sqrt[18]{\frac{\varepsilon}{\lambda}}
\cdot \left|\log\frac{\varepsilon}{\lambda}\right|^{n}   .
\end{eqnarray}

To replace $h_{1/2}$ by $h$  in  \eqref{fh12} and \eqref{gh12}, we claim that
\begin{equation}
\label{hlambdah12}
\int_{\R^n}|h-h_{1/2}| \leq 5\tilde{c} n^8 
\sqrt[18]{\frac{\varepsilon}{\lambda}}
\cdot \left|\log\frac{\varepsilon}{\lambda}\right|^{n}  .
\end{equation}
To prove \eqref{hlambdah12}, we consider
\begin{eqnarray*}
X_-&=&
\{x\in \R^n: h(x)\leq h_{1/2}(x) \}\\
X_+&=&
\{x\in \R^n: h(x)> h_{1/2}(x) \}.
\end{eqnarray*}
It follows from \eqref{hlarger-tildefg} that for any $x\in X_-$, we have
$$
h(x)\geq \tilde{f}(x)^{1-\lambda}\tilde{g}(x)^{\lambda}\geq \min\{\tilde{f}(x),\tilde{g}(x)\},
$$
or in other words, if $x\in X_-$, then
$$
0\leq h_{1/2}(x)-h(x)\leq |h_{1/2}(x)-\tilde{f}(x)|+ |h_{1/2}(x)-\tilde{g}(x)|.
$$
In particular, \eqref{fh12} and \eqref{gh12} imply
\begin{eqnarray}
\nonumber
\int_{X_-}|h-h_{1/2}|&=&
\int_{X_-}(h_{1/2}-h)\leq \int_{X_-}(|h_{1/2}-\tilde{f}|+ |h_{1/2}-\tilde{g}|)\\
\label{hlambdah12sandwichfg}
&\leq &2\tilde{c} n^8 \sqrt[18]{\frac{\varepsilon}{\lambda}}
\cdot \left|\log\frac{\varepsilon}{\lambda}\right|^{n}.
\end{eqnarray}
On the other hand, $\int_{\R^n}h<1+\varepsilon$
and
$\int_{\R^n} h_{1/2}\geq 1$ by \eqref{PLhstabcor-cond2}, thus \eqref{hlambdah12sandwichfg} implies
\begin{eqnarray}
\nonumber
\int_{X_+}|h-h_{1/2}|&=& \int_{X_+}(h-h_{1/2})=
\int_{\R^n}h-\int_{\R^n} h_{1/2}+
\int_{X_-}(h_{1/2}-h)\\
\label{hlambda-larger-h12}
&\leq& \varepsilon+\int_{X_-}(h_{1/2}-h)
\leq 3\tilde{c} n^8 \sqrt[18]{\frac{\varepsilon}{\lambda}}
\cdot \left|\log\frac{\varepsilon}{\lambda}\right|^{n}.
\end{eqnarray}
We conclude \eqref{hlambdah12} by \eqref{hlambdah12sandwichfg} and \eqref{hlambda-larger-h12}.

To replace $\tilde{f}$ and $\tilde{g}$ by $f$ and $g$ in  \eqref{fh12} and \eqref{gh12}, we claim that
\begin{equation}
\label{ftildef-gtildeg}
\int_{\R^n}|f-\tilde{f}|\leq \varepsilon\mbox{ \ and \ }
\int_{\R^n}|g-\tilde{g}|\leq \varepsilon.
\end{equation}
Readily $\tilde{f}\geq f$ and $\tilde{g}\geq g$. 
It follows from \eqref{PLhstabcor-cond2} and $\int_{\R^n}g\leq \int_{\R^n}\tilde{g}=1$ that
$$
\int_{\R^n}|f-\tilde{f}|=\int_{\R^n}\tilde{f}-\int_{\R^n}f\leq 1-\frac1{1+\varepsilon}<\varepsilon,
$$
and similar argument for $g$ and $\tilde{g}$ completes the proof of \eqref{ftildef-gtildeg}.

We conclude from \eqref{fh12}, \eqref{gh12}, \eqref{hlambdah12} and \eqref{ftildef-gtildeg}
that
\begin{eqnarray*}
\int_{\R^n}|f-h|\,dx&\leq &7\tilde{c} n^8\sqrt[18]{\frac{\varepsilon}{\lambda}}
\cdot \left|\log\frac{\varepsilon}{\lambda}\right|^{n} , \\
\int_{\R^n}|g-h|\,dx&\leq &7\tilde{c} n^8\sqrt[18]{\frac{\varepsilon}{\lambda}}
\cdot \left|\log\frac{\varepsilon}{\lambda}\right|^{n}   ,
\end{eqnarray*}
proving Theorem~\ref{PLhstab}.
\proofbox

\section{Proof of Theorem~\ref{PLhstab0} and Corollary~\ref{PLhstablikeMaggi}}
\label{secThCor}

First we verify a simple estimate.

\begin{lemma}
\label{logtn}
If $\varrho>0$, $t>1$ and $n\geq 2$, then
$$
(\log t)^n\leq \left(\frac{n\varrho}{e}\right)^n t^{\frac1{\varrho}}.
$$
\end{lemma}
\proof We observe that for $s=\log t$, differentiating the function $s\mapsto n\log s-\frac{s}{\varrho}$ yields
$$
\log\frac{(\log t)^n}{t^{\frac1{\varrho}}}=n\log s-\frac{s}{\varrho}\leq n(\log (n\varrho)-1)=n\log\frac{n\varrho}{e}.
$$
In turn, we conclude Lemma~\ref{logtn}.
\proofbox

\noindent{\bf Proof of Theorem~\ref{PLhstab0} and Corollary~\ref{PLhstablikeMaggi}: }
We may assume that $f$ and $g$ are probability densities.

For Theorem~\ref{PLhstab0}, we deduce from Theorem~\ref{PLhstab} and Lemma~\ref{logtn} the following statement:
For some absolute constants $c_1,c_2>1$, if  $\varepsilon<c_1^{-n}n^{-n}\cdot \tau$, then
there exists $w\in\R^n$ such that 
\begin{eqnarray}
\label{19Phstabf}
\int_{\R^n}|f(x)- h(x-\lambda\,w)|\,dx&\leq &c_2^nn^{n}\sqrt[19]{\frac{\varepsilon}{\tau}}\\
\label{19Phstabg}
\int_{\R^n}|g(x)-h(x+(1-\lambda)w)|\,dx&\leq &c_2^nn^n\sqrt[19]{\frac{\varepsilon}{\tau}},
\end{eqnarray}
settling Theorem~\ref{PLhstab0} if $\varepsilon<c_1^{-n}n^{-n}\cdot \tau$.

On the other hand,  if  $\varepsilon\geq c_1^{-n}n^{-n}\cdot \tau$, then
we use that the left hand sides of \eqref{19Phstabf} and \eqref{19Phstabg} are
at most $2+\varepsilon\leq 3$; therefore both \eqref{19Phstabf} and \eqref{19Phstabg} readily hold for suitable
absolute constant $c_2>1$. This proves Theorem~\ref{PLhstab0}.

For Corollary~\ref{PLhstablikeMaggi}, the functions
$f$ and $g$ are log-concave probability densities on $\R^n$. In this case, we define
$$
h(z)=\sup_{z=(1-\lambda)x+\lambda y}f(x)^{1-\lambda}g(y)^\lambda,
$$
which is log-concave on $\R^n$ according to Lemma~\ref{log-concave-geometric-mean} (i).
For the $w$ in \eqref{19Phstabf} and \eqref{19Phstabg}, 
we deduce  that
$$
\widetilde{L}_1(f,g)\leq \int_{\R^n}|f(x+w)-g(x)|\,dx\leq 2c_2^nn^n\sqrt[19]{\frac{\varepsilon}{\tau}},
$$
yielding Corollary~\ref{PLhstablikeMaggi}.
\proofbox

\section{Proof of Theorem~\ref{PLhstabcor0}}
\label{secmanyfunctions}

We deduce from Lemma~\ref{log-concave-geometric-mean} (i) and induction on $m$ the following corollary.

\begin{coro}
\label{log-concave-geometric-mean-cor}
If $\lambda_1,\ldots,\lambda_m>0$ satisfy $\sum_{i=1}^m\lambda_i=1$ and 
$f_1,\ldots,f_m$ are log-concave functions with positive integral on $\R^n$, then
$$
h(z)=\sup_{z=\sum_{i=1}^m\lambda_ix_i} \prod_{i=1}^mf(x_i)^{\lambda_i}
$$
is log-concave and has positive integral.
\end{coro}

The first main step towards proving Theorem~\ref{PLhstabcor} is the case when each $\lambda_i$ in Theorem~\ref{PLhstabcor} is $\frac1m$.

\begin{theo}
\label{PLhstabcor1m}
Let $c>1$ be the absolute constant in Theorem~\ref{PLhstab}, let $\gamma_0=cn^8$ and 
$\varepsilon_0=c^{-n}n^{-n}$.
If  
$f_1,\ldots,f_m$, $m\geq 2$ are log-concave probability densities on $\R^n$ such that
$$
\int_{\R^n}\sup_{mz=\sum_{i=1}^m x_i} \prod_{i=1}^mf_i(x_i)^{\frac1m}\,dz
\leq 1+\varepsilon  
$$
for $0<\varepsilon<\varepsilon_0/m^4$, then for the log-concave 
$h(z)=\sup_{mz=\sum_{i=1}^mx_i} \prod_{i=1}^mf(x_i)^{\frac1m}$,
there exist  $w_1,\ldots,w_m\in\R^n$ such that  $\sum_{i=1}^mw_i=o$ and
$$
\int_{\R^n}|f_i(x)-h(x+w_i)|\,dx\leq m^4\cdot\gamma_0
\sqrt[18]{\varepsilon}\cdot
\left|\log\varepsilon\right|^{n}.
$$
\end{theo}
\proof Instead of Theorem~\ref{PLhstabcor1m},
we prove that if  $0<\varepsilon<\varepsilon_0/4^{\lceil\log_2m\rceil}$, then 
there exist  $w_1,\ldots,w_m\in\R^n$ such that  $\sum_{i=1}^mw_i=o$ and
\begin{equation}
\label{PLhstabcor1m0}
\int_{\R^n}|f_i(x)-h(x+w_i)|\,dx\leq 4^{\lceil\log_2m\rceil}\cdot\gamma_0
\sqrt[18]{\varepsilon}\cdot
\left|\log\varepsilon\right|^{n} .
\end{equation}
Since $4^{\lceil\log_2m\rceil}<4^{2\log_2m}=m^4$, \eqref{PLhstabcor1m0} yields 
Theorem~\ref{PLhstabcor1m}.

We prove \eqref{PLhstabcor1m0}  by induction on $\lceil\log_2m\rceil\geq 1$. 
If 
$\lceil\log_2m\rceil=1$, and hence $m=2$, then 
$$
h(z)=\sup_{z=\lambda_1x_1+\lambda_2x_2} f_1(x_1)^{\lambda_1}f_2(x_2)^{\lambda_2}
$$
is log-concave by Lemma~\ref{log-concave-geometric-mean} (i). Therefore, the case $m=2$ of
\eqref{PLhstabcor1m0} follows from Theorem~\ref{PLhstab}.

Next, we assume that $\lceil\log_2m\rceil>1$, and let $k=\lceil m/2\rceil$, and hence 
$\lceil\log_2(m-k)\rceil\leq\lceil\log_2k\rceil=\lceil\log_2m\rceil-1$. 
We consider the coefficient
$$
\lambda=\frac{m-k}m \mbox{ \ satisfying \ }\frac13\leq\lambda\leq \frac12,
$$ 
and the functions
\begin{eqnarray}
\label{cor15hdef}
h(z)&=&\sup_{mz=\sum_{i=1}^mx_i} \prod_{i=1}^mf_i(x_i)^{\frac1m}\\
\label{cor15gdef}
f(z)&=&\sup_{kz=\sum_{i=1}^kx_i} \prod_{i=1}^kf_i(x_i)^{\frac1k}\\
\label{cor15fdef}
g(z)&=&\sup_{(m-k)z=\sum_{i=k+1}^mx_i} \prod_{i=k+1}^mf_i(x_i)^{\frac1{m-k}},
\end{eqnarray}
which are log-concave by Corollary~\ref{log-concave-geometric-mean-cor}.
In particular, we have
$$
h(z)=\sup_{z=\lambda x+(1-\lambda)y} f(x)^{1-\lambda}g(y)^{\lambda}.
$$
It follows from the Pr\'ekopa-Leindler inequality that
\begin{eqnarray}
\label{cor15flow}
\int_{\R^n}f&\geq &1\\
\label{cor15glow}
\int_{\R^n}g&\geq &1\\
\label{cor15h}
\int_{\R^n}h&\geq &\left(\int_{\R^n}f\right)^{1-\lambda}\left(\int_{\R^n}g\right)^{\lambda}\geq 1.\\
\end{eqnarray}
Since $\int_{\R^n}h<1+\varepsilon$, we deduce on the one hand, that
\begin{eqnarray}
\label{cor15fupp}
\int_{\R^n}f&\leq & (1+\varepsilon)^{\frac1{1-\lambda}}\leq (1+\varepsilon)^3\leq 1+4\varepsilon\\
\label{cor15gupp}
\int_{\R^n}g&\leq &  1+4\varepsilon,
\end{eqnarray}
and on the other hand, Theorem~\ref{PLhstab} yields  that 
for $a=\int_{\R^n}g/\int_{\R^n}f$, 
there exists $w\in\R^n$ such that
\begin{eqnarray}
\nonumber
\int_{\R^n}|f(x)-a^{\lambda}h(x-\lambda\,w)|\,dx&\leq &\gamma_0\sqrt[18]{\frac{\varepsilon}{1/3}}\cdot|\log\varepsilon|^{n}
\int_{\R^n} f\\
\label{afhcor15}
&\leq& 2\gamma_0\sqrt[18]{\varepsilon}\cdot|\log\varepsilon|^{n}\\
\nonumber
\int_{\R^n}|g(x)-a^{-(1-\lambda)}h(x+(1-\lambda)w)|\,dx&\leq &\gamma_0\sqrt[18]{\frac{\varepsilon}{1/3}}\cdot|\log\varepsilon|^{n} \int_{\R^n} g\\
\label{aghcor15}
&\leq& 2 \gamma_0\sqrt[18]{\varepsilon}\cdot|\log\varepsilon|^{n}.
\end{eqnarray}
We deduce from \eqref{cor15flow}, \eqref{cor15glow}, \eqref{cor15fupp}, \eqref{cor15gupp} that
$1+4\varepsilon>a,a^{-1}>\frac1{1+4\varepsilon}>1-4\varepsilon$; therefore,
$\frac13\leq\lambda\leq\frac23$,
$\int_{\R^n}h<1+\varepsilon$, \eqref{afhcor15} and \eqref{aghcor15} yield
\begin{eqnarray}
\label{fhcor15}
\int_{\R^n}|f(x)-h(x-\lambda\,w)|\,dx&\leq &4\gamma_0\sqrt[18]{\varepsilon}\cdot|\log\varepsilon|^{n}\\
\label{ghcor15}
\int_{\R^n}|g(x)-h(x+(1-\lambda)w)|\,dx&\leq &4\gamma_0\sqrt[18]{\varepsilon}\cdot|\log\varepsilon|^{n} .
\end{eqnarray}

Since $\lceil\log_2(m-k)\rceil\leq\lceil\log_2k\rceil=\lceil\log_2m\rceil-1$, induction
and \eqref{cor15fdef}, \eqref{cor15gdef} and \eqref{cor15hdef} yield that there exist
$\tilde{w}_1,\ldots,\tilde{w}_m\in\R^n$ such that  
\begin{equation}
\label{tildewsum}
\sum_{i=1}^k\tilde{w_i}=\sum_{j=k+1}^m\tilde{w_j}=o,
\end{equation}
and if $i=1,\ldots,k$ and $j=k+1,\ldots,m$, then
\begin{eqnarray}
\label{fifcor15}
\int_{\R^n}|f_i(x)-f(x+\tilde{w}_i)|\,dx&\leq &4^{\lceil\log_2 k\rceil}\gamma_0\sqrt[18]{4\varepsilon}\cdot|\log\varepsilon|^{n},\\
\label{fjgcor15}
\int_{\R^n}|f_j(x)-g(x+\tilde{w}_j)|\,dx&\leq &4^{\lceil\log_2 (m-k)\rceil}\gamma_0\sqrt[18]{4\varepsilon}\cdot|\log\varepsilon|^{n}.
\end{eqnarray}
Combining \eqref{fhcor15}, \eqref{ghcor15}, \eqref{fifcor15} and  \eqref{fjgcor15}
shows the existence of $w_1,\ldots,w_m\in\R^n$ such that if $i=1,\ldots,m$, then 
$$
\int_{\R^n}|f_i(x)-h(x+w_i)|\,dx\leq 
4\cdot 4^{\lceil\log_2 m\rceil-1}\gamma_0\sqrt[18]{\varepsilon}\cdot|\log\varepsilon|^{n}=4^{\lceil\log_2 m\rceil}
 \gamma_0\sqrt[18]{\varepsilon}\cdot|\log\varepsilon|^{n};
$$
 namely,
\begin{eqnarray*}
w_i&= &-\lambda\,w-\tilde{w}_i\mbox{ \  for $i=1,\ldots,k$},\\
w_j&= &(1-\lambda)w-\tilde{w}_j\mbox{ \  for $j=k+1,\ldots,m$}.
\end{eqnarray*}
Since $\lambda=\frac{m-k}m$ and $1-\lambda=\frac{k}m$, we have
$$
\sum_{i=1}^mw_i=-k\cdot \lambda\,w-\left(\sum_{i=1}^k\tilde{w_i}\right)
+(m-k)(1-\lambda)w-\left(\sum_{j=k+1}^m\tilde{w_j}\right)=o,
$$
proving \eqref{PLhstabcor1m0}.
\proofbox

For $m\geq 2$, we consider the $(m-1)$-simplex
$$
\Delta^{m-1}=\{p=(p_1,\ldots,p_m)\in\R^m:\,p_1+\ldots+p_m=1\}.
$$
The proof of Lemma~\ref{log-concave-geometric-mean} readily extends to verify he following statement.

\begin{lemma}
\label{log-concave-geometric-meanm}
Let $f_1,\ldots,f_m$, $m\geq 2$ be log-concave probability densities  with positive integral on $\R^n$, $n\geq 2$.
For $p=(p_1,\ldots,p_m)\in\Delta^{m-1}$, the function
$$
h_p(z)=\sup_{z=\sum_{i=1}^m p_ix_i} \prod_{i=1}^mf_i(x_i)^{p_i}
$$
is log-concave on $\R^n$, and 
 the function $p\mapsto\int_{\R^n}h_p$
is log-concave on $\Delta^{m-1}$.
\end{lemma}

\begin{theo}
\label{PLhstabcor}
For some absolute constant $\tilde{\gamma}>1$, if  $\tau\in(0,\frac1m]$, $m\geq 2$, 
$\lambda_1,\ldots,\lambda_m\in[\tau,1-\tau]$ satisfy $\sum_{i=1}^m\lambda_i=1$ and 
$f_1,\ldots,f_m$ are log-concave functions with positive integral on $\R^n$ such that
$$
\int_{\R^n}\sup_{z=\sum_{i=1}^m\lambda_ix_i} \prod_{i=1}^mf_i(x_i)^{\lambda_i}\,dz
\leq (1+\varepsilon)  \prod_{i=1}^m \left(\int_{\R^n}f_i\right)^{\lambda_i}
$$
for $0<\varepsilon<\tau\cdot \tilde{\gamma}^{-n}n^{-n}/m^4$, then for the log-concave 
$h(z)=\sup_{z=\sum_{i=1}^m\lambda_ix_i} \prod_{i=1}^mf(x_i)^{\lambda_i}$,
there exist $a_1,\ldots,a_m>0$ and  $w_1,\ldots,w_m\in\R^n$ such that
$\sum_{i=1}^m\lambda_iw_i=o$ and
for $i=1,\ldots,m$, we have
$$
\int_{\R^n}|f_i(x)-a_ih(x+w_i)|\,dx\leq \tilde{\gamma}m^5n^8
\sqrt[18]{\frac{\varepsilon}{m\tau}}\cdot
\left|\log\frac{\varepsilon}{m\tau}\right|^{n} 
\int_{\R^n}f_i.
$$
\end{theo}
\noindent{\bf Remark } $a_i=\frac{\left(\int_{\R^n}f_i\right)^{1-\lambda_i}}
{\prod_{j\neq i}\left(\int_{\R^n}f_j\right)^{\lambda_j}}$ for $i=1,\ldots,m$
in Theorem~\ref{PLhstabcor}.\\
\proof
Let $\tau\in(0,\frac1m]$, let $\lambda_1,\ldots,\lambda_m\in[\tau,1-\tau]$ with $\lambda_1+\ldots+\lambda_m=1$ and let $f_1,\ldots,f_m$ be
log-concave with positive integralral as in Theorem~\ref{PLhstabcor}. In particular,
$\lambda=(\lambda_1,\ldots,\lambda_m)\in \Delta^{m-1}$. We may assume that 
$f_1,\ldots,f_m$ are probability densities, $\lambda_1\geq \ldots\geq \lambda_m$ and $\lambda_m<\frac1m$ 
(as if $\lambda_m=\frac1m$, then Theorem~\ref{PLhstabcor1m} implies Theorem~\ref{PLhstabcor}).

Let $\tilde{p}=(\frac1m,\ldots,\frac1m)$, and for
$i=1,\ldots,m$, let $v_{(i)}\in\R^m$ be the vector whose $i$th coordinate is $1$, and the rest is $0$, and hence 
$v_{(1)},\ldots,v_{(m)}$ are the vertices of $\Delta^{m-1}$, and
$\lambda=\sum_{i=1}^m\lambda_iv_{(i)}$. For $p=(p_1,\ldots,p_m)\in\Delta^{m-1}$, we write
$$
h_p(z)=\sup_{z=\sum_{i=1}^m p_ix_i} \prod_{i=1}^mf_i(x_i)^{p_i},
$$
and hence $h_{v_{(i)}}=f_i$. According the conditions in Theorem~\ref{PLhstabcor},
\begin{equation}
\label{hlambdaepscond}
\int_{\R^n}h_\lambda<1+\varepsilon.
\end{equation}

Since $\lambda_i-\lambda_m\geq 0$ for $i=1,\ldots,m-1$ and $m\lambda_m<1$, it follows that
$$
q=\sum_{i=1}^{m-1}\frac{\lambda_i-\lambda_m}{1-m\lambda_m}\cdot v_{(i)}\in \Delta^{m-1},
$$
and Corollary~\ref{log-concave-geometric-meanm} and $\int_{\R^n}h_{v_{(i)}}=\int_{\R^n}f_i=1$
 yield that $\int_{\R^n}h_q\geq 1$. Since
$$
\lambda=m\lambda_m\tilde{p}+\sum_{i=1}^{m-1}(\lambda_i-\lambda_m)v_{(i)}=
m\lambda_m\tilde{p}+(1-m\lambda_m)q,
$$ 
\eqref{hlambdaepscond} and 
Corollary~\ref{log-concave-geometric-meanm} imply that
$$
1+\varepsilon> \int_{\R^n}h_\lambda\geq \left(\int_{\R^n}h_q\right)^{1-m\lambda_m}
\left(\int_{\R^n}h_{\tilde{p}}\right)^{m\lambda_m}\geq 
\left(\int_{\R^n}h_{\tilde{p}}\right)^{m\lambda_m};
$$
therefore, $\varepsilon\leq m\tau\leq m\lambda_m$ yields
$$
\int_{\R^n}h_{\tilde{p}}<(1+\varepsilon)^{\frac1{m\lambda_m}}<
e^{\frac{\varepsilon}{m\lambda_m}}<1+\frac{2\varepsilon}{m\lambda_m}.
$$
According to Theorem~\ref{PLhstabcor1m}, there exist  $w_1,\ldots,w_m\in\R^n$ such that  $\sum_{i=1}^mw_i=o$ and
$$
\int_{\R^n}|f_i(x+w_i)-h_{\tilde{p}}(x)|\,dx\leq m^4\cdot\gamma_0
\sqrt[18]{\frac{2\varepsilon}{m\lambda_m}}\cdot
\left|\log\frac{2\varepsilon}{m\lambda_m}\right|^{n}
$$
for $i=1,\ldots,m$. Replacing $f_i(x)$ by $f_i(x+w_i)$ for $i=1,\ldots,m$ does not change $h_{\tilde{p}}$ by the condition
$\sum_{i=1}^mw_i=o$; therefore, we may assume that
\begin{equation}
\label{fihtildep}
\int_{\R^n}|f_i(x)-h_{\tilde{p}}(x)|\,dx\leq m^4\cdot\gamma_0
\sqrt[18]{\frac{2\varepsilon}{m\lambda_m}}\cdot \left|\log\frac{2\varepsilon}{m\lambda_m}\right|^{n}
\end{equation}
for $i=1,\ldots,m$.

To replace $h_{\tilde{p}}$ by $h_\lambda$  in  \eqref{fihtildep}, we claim that
\begin{equation}
\label{hlambdahtildep}
\int_{\R^n}|h_\lambda-h_{\tilde{p}}| \leq 
3 m^5\cdot\gamma_0
\sqrt[18]{\frac{2\varepsilon}{m\lambda_m}}\cdot
\left|\log\frac{2\varepsilon}{m\lambda_m}\right|^{n}  .
\end{equation}
To prove \eqref{hlambdahtildep}, we consider
\begin{eqnarray*}
X_-&=&
\{x\in \R^n: h_\lambda(x)\leq h_{\tilde{p}}(x) \}\\
X_+&=&
\{x\in \R^n: h_\lambda(x)> h_{\tilde{p}}(x) \}.
\end{eqnarray*}
It follows from the definition of $h_\lambda$ that for any $x\in X_-$, we have
$$
h_\lambda(x)\geq \prod_{i=1}^mf_i(x)^{\lambda_i}\geq \min\{f_1(x),\ldots,f_m(x)\},
$$
or in other words, if $x\in X_-$, then
$$
0\leq h_{\tilde{p}}(x)-h_\lambda(x)\leq  \sum_{i=1}^m|f_i(x)-h_{\tilde{p}}(x)|.
$$
In particular, \eqref{fihtildep} implies
\begin{eqnarray}
\nonumber
\int_{X_-}|h_\lambda-h_{\tilde{p}}|&=&
\int_{X_-}(h_{\tilde{p}}-h_\lambda)\leq \sum_{i=1}^m\int_{X_-}|f_i(x)-h_{\tilde{p}}(x)|\\
\label{hlambdahtildepsandwich}
&\leq & m^5\cdot\gamma_0
\sqrt[18]{\frac{2\varepsilon}{m\lambda_m}}\cdot \left|\log\frac{2\varepsilon}{m\lambda_m}\right|^{n}.
\end{eqnarray}

On the other hand, $\int_{\R^n}h_\lambda<1+\varepsilon$
and the Pr\'ekopa-Leindler inequality yields
$\int_{\R^n} h_{\tilde{p}}\geq 1$, thus \eqref{hlambdahtildepsandwich} implies
\begin{eqnarray}
\nonumber
\int_{X_+}|h_\lambda-h_{\tilde{p}}|&=& \int_{X_+}(h_\lambda-h_{\tilde{p}})=
\int_{\R^n}h_\lambda-\int_{\R^n} h_{\tilde{p}}+
\int_{X_-}(h_{\tilde{p}}-h_\lambda)\\
\label{hlambda-larger-htildep}
&\leq& \varepsilon+\int_{X_-}(h_{\tilde{p}}-h_\lambda)
\leq  2m^5\cdot\gamma_0
\sqrt[18]{\frac{2\varepsilon}{m\lambda_m}}\cdot \left|\log\frac{2\varepsilon}{m\lambda_m}\right|^{n}.
\end{eqnarray}
We conclude \eqref{hlambdahtildep} by \eqref{hlambdahtildepsandwich} and \eqref{hlambda-larger-htildep}.

Finally, combining \eqref{fihtildep} and \eqref{hlambdahtildep} prove Theorem~\ref{PLhstabcor}.
\proofbox

\noindent{\bf Proof of Theorem~\ref{PLhstabcor0} } 
We may assume that
$\int_{\R^n}f_i=1$ for $i=1,\ldots,m$ in Theorem~\ref{PLhstabcor0}
for the log-concave functions $f_1,\ldots,f_m$ on $\R^n$.

Let  $\tau\in(0,\frac1m]$ for $m\geq 2$, and let
$\lambda_1,\ldots,\lambda_m\in[\tau,1-\tau]$ satisfy $\sum_{i=1}^m\lambda_i=1$ 
 such that
$$
\int_{\R^n}\sup_{z=\sum_{i=1}^m\lambda_ix_i} \prod_{i=1}^mf_i(x_i)^{\lambda_i}\,dz
\leq 1+\varepsilon
$$
for $\varepsilon\in(0,1]$.

For the absolute constant $\tilde{\gamma}>1$ of Theorem~\ref{PLhstabcor},
if 
\begin{equation}
\label{PLhstabcoreps}
0<\varepsilon<\tau\cdot \tilde{\gamma}^{-n}n^{-n}/m^4, 
\end{equation}
then for the log-concave 
$h(z)=\sup_{z=\sum_{i=1}^m\lambda_ix_i} \prod_{i=1}^mf(x_i)^{\lambda_i}$,
there exist   $w_1,\ldots,w_m\in\R^n$ such that
$\sum_{i=1}^m\lambda_iw_i=o$ and
for $i=1,\ldots,m$, we have
$$
\int_{\R^n}|f_i(x)-h(x+w_i)|\,dx\leq \tilde{\gamma}m^5n^8
\sqrt[18]{\frac{\varepsilon}{m\tau}}\cdot
\left|\log\frac{\varepsilon}{m\tau}\right|^{n}.
$$

We deduce from Lemma~\ref{logtn} that
\begin{equation}
\label{PLhstabcorest}
\int_{\R^n}|f_i(x)-h(x+w_i)|\,dx\leq \tilde{\gamma}_0^nn^nm^5
\sqrt[19]{\frac{\varepsilon}{m\tau}}
\end{equation}
for $i=1,\ldots,m$ and some absolute constant $\tilde{\gamma}_0\geq \max\{\tilde{\gamma},3\}$, proving Theorem~\ref{PLhstabcor0} if \eqref{PLhstabcoreps} holds. Finally
if $\varepsilon\geq \tau\cdot \tilde{\gamma}^{-n}n^{-n}/m^4$, then 
\eqref{PLhstabcorest} readily holds as the left hand side is at most $2+\varepsilon\leq 3$.
\proofbox

\noindent{\bf Acknowledgement} 
We are grateful 
for  useful discussions with Galyna Livshyts, Alessio Figalli, Emanuel Milman, Ronen Eldan and Marco Barchiesi.

\noindent K\'aroly J. B\"or\"oczky\\
Alfr\'ed R\'enyi Institute of Mathematics, Re\'altanoda u. 13-15, H-1053 Budapest, Hungary, and\\
Department of Mathematics, Central European University, N\'ador u. 9, H-1051, Budapest, Hungary\\
boroczky.karoly.j@renyi.hu\\[4ex]
Apratim De\\
Department of Mathematics, Central European University, N\'ador u. 9, H-1051, Budapest, Hungary\\
de.apratim91@gmail.com

\end{document}